\documentclass[10pt,A4]{amsart}
\usepackage{amssymb,amsmath,amsthm}
\usepackage[all]{xy}
\usepackage{epsfig,exscale}
\usepackage{xspace}
\usepackage{color}
\usepackage{colordvi}
\usepackage[small,nohug,heads=littlevee]{diagrams}
\newcommand{\nc}{\newcommand}

\nc{\id}{{\mathrm{id}}}

\def\sha{{\mbox{\cyr X}}}

\newtheorem{thm}{Theorem}
\newtheorem{cor}[thm]{Corollary}

\newtheorem{lem}[thm]{Lemma}
\newtheorem{prop}[thm]{Proposition}
\newtheorem{defn}{Definition}
\newtheorem{rmk}[thm]{Remark}
\newtheorem{example}[thm]{Example}

\font \eightrm=cmr8

\font\cyr=wncyr10

\begin{document}

\title{Generalized Shuffles Related to Nijenhuis and $TD$-Algebras}

\author{Kurusch Ebrahimi-Fard}
\address{I.H.\'E.S.,
         Le Bois-Marie,
         35, Route de Chartres,
         F-91440 Bures-sur-Yvette, France}
         \email{kurusch@ihes.fr}
         \urladdr{http://www.th.physik.uni-bonn.de/th/People/fard/}

\author{Philippe Leroux}
\address{Anciennement rattach\'e \`a l'I.R.M.A.R., Universit\'e de Rennes
         I et C.N.R.S. U.M.R. 6625 campus de Beaulieu, 35042 Rennes Cedex France. }
         \email{ ph$\_$ler$\_$math@yahoo.com }

\date{August 19, 2008}

\maketitle

\begin{abstract}
Shuffle type products are well-known in mathematics and
physics. They are intimately related to Loday's dendriform
algebras and were extensively used to give explicit constructions
of free Rota--Baxter algebras. In the literature there exist at
least two other Rota--Baxter type algebras, namely, the Nijenhuis
algebra and the so-called $TD$-algebra. The explicit construction
of the free unital commutative Nijenhuis algebra uses a modified
quasi-shuffle product, called the right-shift shuffle. We show
that another modification of the quasi-shuffle, the so-called
left-shift shuffle, can be used to give an explicit construction
of the free unital commutative $TD$-algebra. We explore some basic
properties of $TD$-operators. Our construction is related to
Loday's unital commutative tridendriform algebra, including the
involutive case. The concept of Rota--Baxter, Nijenhuis and
$TD$-bialgebras is introduced at the end and we show that any
commutative bialgebra provides such objects.
\end{abstract}


\section{Introduction}
\label{sect:intro}

Through the connection between free commutative Rota--Baxter
algebras and Hoffman's quasi-shuffle algebra~\cite{EGmixShuf}, the
development of further instances of generalized shuffle products
for other algebras which are characterized by a particular
operator identity, e.g. of Rota--Baxter type, is an interesting
exercise. It provides further insights into such algebras as well
as their relations to other fields in mathematics and physics, see
e.g.~\cite{KreimerShuf,Rosso}. For instance, from the link between
Hoffman's quasi-shuffle product~\cite{Hoffman00} and the mixable
shuffle product~\cite{GK00} one realizes that Hoffman's well-known
partition identity for multiple-zeta-values~\cite{Hoffman05} is
another instance of Spitzer's classical identity for commutative
Rota--Baxter algebras~\cite{EGmzv}.

In the recent literature further Rota--Baxter type algebras
appeared, to wit, the associative Nijenhuis algebra and the
so-called $TD$-algebra. In fact, the latter is a particular
Rota--Baxter algebra of generalized weight. Recall that the Lie
algebraic version of the Nijenhuis as well as Rota--Baxter
relation are well-known in the theory of classical integrable
systems~\cite{BelavinDrinfeld82,GolSok00,KSM90,STS83}. However,
$TD$-algebras only recently entered the scene in the work of one
of us~\cite{Leroux04b} in the context of Loday's dendriform
algebras~\cite{Loday01}. In this work we will show how
$TD$-algebras relate to proper Rota--Baxter algebras of
generalized, that is, non-scalar weight.

The construction of free commutative Rota--Baxter algebras
appeared first in the work of Baxter~\cite{Baxter60} and
Rota~\cite{Rota69}. Later, Cartier~\cite{Cartier72} gave a
description which already contained the notion of quasi-shuffle
product. Guo and Keigher~\cite{GK00} established this link most
clearly by following closely the work of the former authors.
Following~\cite{GK00} one of us introduced a modified
quasi-shuffle product, called right-shift shuffle, to construct
the free commutative Nijenhuis algebra~\cite{KEF2}. In this paper
we continue this line of thought and establish the free
commutative $TD$-algebra, using another modified quasi-shuffle
product, called left-shift shuffle.

Let us briefly outline the organization of the paper. In
Section~\ref{sect:dendal} we remind the reader of the definitions
of dendriform algebras. Section~\ref{sect:RBtypeA} summarizes
briefly the theory of Rota--Baxter algebras of non-scalar weight
as well as that of Nijenhuis algebras. The link between
$TD$-algebras and Rota--Baxter algebras of generalized weight is
established. The categories of Rota--Baxter and Nijenhuis algebras
are described in Subsection~\ref{ssect:RBcategories}. In
Section~\ref{sect:shuffle} we recall for completeness three
different types of shuffle products existing in the literature and
introduce the so-called left-shift shuffle. These shuffles are the
main ingredients for the construction of free objects in the
categories defined in Subsection~\ref{ssect:RBcategories}. We show
that the free commutative Nijenhuis algebra is a $TD$-algebra and
construct the free commutative $TD$-algebra. Relations with
commutative tridendriform algebras are exhibited including the
involutive cases. In Section~\ref{sect:QshufLshuf} we present a
relation between the quasi-shuffle product and the left-shift
shuffle product. This allows us to present another proof of a
result of Loday on the free commutative tridendriform algebra. The
next section contains the notion of Rota--Baxter type bialgebra.
We construct a covariant functor from the category of commutative
bialgebras to the category of Rota--Baxter type bialgebras thanks
to the free objects obtained earlier and show that this result
survives a priori only for the Nijenhuis case if we start with the
category of commutative algebras instead.

\smallskip

In the sequel, $\mathbb{K}$, with $char(\mathbb{K})=0$, denotes
the ground field over which all algebraic structures are defined.


\section{Dendriform algebras}
\label{sect:dendal}

Dendriform algebras were introduced by Loday~\cite{Loday01} with
motivations coming from K-Theory. It was shown that the free
dendriform algebra on one generator can be constructed over planar
binary rooted trees. To extend this construction to planar rooted
trees, Loday and Ronco \cite{LodayRonco04} introduced
tridendriform algebras.

\begin{defn} \cite{LodayRonco04} \label{def:dendtri}
A {\bf{tridendriform algebra}} is a $\mathbb{K}$-vector space $T$
equipped with three binary operations, $\prec~$,  $\succ$, and
$\bullet$ from $T^{\otimes 2} \xrightarrow{} T$ satisfying the
following set of relations for all $x,y,z \in T$
$$
    (x \prec y )\prec z = x \prec(y \star z), \
    (x \succ y )\prec z = x \succ(y \prec z), \
    (x \star y )\succ z = x \succ(y \succ z),
$$
$$
    (x \succ y )  \bullet z = x \succ (y \bullet z),  \
    (x \prec y )  \bullet z = x \bullet (y \succ z),  \
    (x \bullet y) \prec z   = x \bullet (y \prec z),  \
    (x \bullet y )\bullet z = x \bullet (y \bullet z),
$$
where by definition
\begin{equation}
\label{def:double}
     x \star y := x  \prec y + x \succ y + x \bullet y,
\end{equation}
for all $x,y \in T$. A tridendriform algebra $T$ is said to be
{\bf{commutative}} if $x \prec y = y \succ x$ and $x \bullet y = y
\bullet x$.
\end{defn}

\noindent The product $\star$ turns out to be associative. If the
composition $\bullet$ in the above definition vanishes or if
$\bullet$ and $\succ$ are combined into the new operation
$\succ':= \succ +\ \bullet$, i.e., $(T,\succ,\prec,\bullet) \to
(T,\succ',\prec)$, then one obtains a {\it{dendriform
algebra}}~\cite{Loday01}. The new product $x \star y := x  \prec y
+ x \succ' y$ is associative as well. For more details, e.g. see
\cite{Aguiar00,Aguiar04,AguLod04,KEF1,EG04,EG05,Leroux04a,Leroux04b,Leroux05}.

Let us denote by \textbf{Dend}, \textbf{TriDend} and
\textbf{CTDend} the category of dendriform, tridendriform and
commutative tridendriform algebras, respectively. We have seen
that there is a canonical functor from the category of
tridendriform algebras to the category of dendriform algebras over
$\mathbb{K}$.

If \textbf{As} and \textbf{Com} denote the category of associative
and commutative algebras, respectively, then we have canonical
functors \textbf{Dend},
\textbf{TriDend}~$\longrightarrow$~\textbf{As} and
\textbf{CTDend}~$\longrightarrow$~\textbf{Com}.


\section{Rota--Baxter type algebras}
\label{sect:RBtypeA}

We now introduce two types of unital algebras, called Rota--Baxter
and Nijenhuis algebra. They are characterized by specific linear
maps each satisfying a particular operator identity. A particular
subclass of Rota--Baxter algebras is identified with
$TD$-algebras. The explicit construction of the free object in the
corresponding categories of commutative algebras (c.f.
Section~\ref{ssect:RBcategories}) denoted by \textbf{ComRBA},
\textbf{ComNA}, respectively, is achieved in terms of generalized
shuffle products.\medskip

In this paper the term algebra always means associative
$\mathbb{K}$-algebra with unit $1_A$, unless otherwise stated. We
denote the product of an algebra $A$ simply by concatenation,
$m_A(a \otimes b)=:ab$. Sometimes we use the bracket notation
i.e., $m_A(a \otimes b)=:[a;b]$.\\


\subsection{Rota--Baxter algebra}
\label{ssect:RBA}

The Rota--Baxter relation was introduced by G.~Baxter in
\cite{Baxter60} motivated by the work of
F.~Spitzer~\cite{Spitzer56}. It was further studied by several
people, among others G.-C.~Rota~\cite{Rota69} and P.~Cartier
\cite{Cartier72}.

\begin{defn} \label{def:RB} Fix $\theta \in
A$. We call the pair $(A,R)$ where $A$ is an algebra and $R: A \to
A$ is a $\mathbb{K}$-linear operator satisfying the
{\bf{Rota--Baxter relation}}
\begin{equation}
    R(x)R(y) = R\bigl(R(x)y + xR(y)\bigr) + R(x \ \theta\ y),\qquad \forall x,y \in A,
    \label{eq:thetaRBR}
\end{equation}
a {\bf{Rota--Baxter algebra}} of weight $\theta$.
\end{defn}
Notice our more general assumption of the weight $\theta$ to be an
element in $A$ in contrast to the definition used in most of the
recent works which restrict $\theta$ to be a scalar parameter in
the underlying field $\mathbb{K}$, see e.g.
\cite{AguLod04,Fields,Leroux04a,Rota95,Rota98,RotaSmith72}.


\subsection{Nijenhuis algebra}
\label{ssect:NA}

In general, the Nijenhuis relation may be regarded as a
homogeneous version of the Rota--Baxter relation of unit weight
$\theta := - 1_A$. In \cite{GolSok00,KSM90} the Lie algebraic
version of the Nijenhuis relation is investigated in the context
of classical Yang--Baxter\footnote{Referring to the physicists
C.~N.~Yang from China and the Australian R.~Baxter.} type
equations, which are also closely related to the Lie algebraic
version of the Rota--Baxter relation of unit weight.

For the Nijenhuis relation in the context of associative algebras,
many of its algebraic aspects were established by Cari{\~n}ena
{\emph{et al.}} in \cite{CGM00}, see also Fuchssteiner's work
\cite{Fuchs97}. In \cite{Leroux04b} a construction of Nijenhuis
operators was proposed from a bialgebraic point of view similar to
a construction for Rota--Baxter operators of weight zero given by
Aguiar, and Aguiar and Loday~\cite{Aguiar04,AguLod04}.

\begin{defn} \label{def:NA}
A {\bf{Nijenhuis algebra}} is a pair $(A,N)$, where $A$ is an
algebra and $N: A \to A$ is a {\bf{Nijenhuis operator}}, that is,
a $\mathbb{K}$-linear operator satisfying the {\bf{Nijenhuis
relation}}
\begin{equation}
    N(x)N(y) = N\bigl(N(x)y + xN(y)\bigr) - N^2(xy),\qquad \forall x,y \in A.
    \label{eq:ANR}
\end{equation}
\end{defn}

Examples of operators fulfilling the Nijenhuis relation may be
found in the cited literature, see especially \cite{CGM00,Leroux04b}.\\


\subsection{General properties of Rota--Baxter type algebras}
\label{ssect:properties}

Let us summarize the main properties of the above-introduced
Rota--Baxter type algebras.

First observe that for $(A,N)$ being a Nijenhuis algebra,
$\tilde{N}:= \id_A - N$ is a Nijenhuis map, too. Both the image of
$N$ and $\tilde{N}$ are subalgebras in $(A,N)$.

In the case of Rota--Baxter algebras we find the following. The
image of a Rota--Baxter map $R$ of weight $\theta \in A$ is a
subalgebra in $A$. One verifies easily the following statement.

\begin{prop} \label{prop:RBdual}
Let $(A,R)$ be a Rota--Baxter algebra of weight $\theta \in A$.
For an element $\theta$ in the center of $A$ $\tilde{R}:= -\theta
\id_A - R$ is a Rota--Baxter map of weight $\theta$ on $A$. The
image of $\tilde{R}$ is a subalgebra in $A$.
\end{prop}


\subsubsection{Double Rota--Baxter type products}
\label{sssect:doubleRBtype}

In the case of a Nijenhuis algebra we recall a result from
\cite{CGM00}. Let $(A,N)$ be a Nijenhuis algebra. We define the
following operation on the vector space underlying $A$
$$
    a \ast_N b := N(a)b + aN(b) - N(ab), \qquad a,b \in A,
$$
which defines a Nijenhuis algebra, denoted by $A_N$, with $N$ as
Nijenhuis map. The pair $(A_N, N)$ is called {\bf{double Nijenhuis
algebra}}.

From the construction and relation~(\ref{eq:ANR}) it is readily
verified that
\begin{equation}
    \label{eq:NAhom}
    N(a \ast_N b) = N(a) N(b)
    \quad {\mathrm{ and }} \quad
    \tilde{N}(a \ast_N b) = - \tilde{N}(a) \tilde{N}(b).
\end{equation}

Let $(A,R)$ be a Rota--Baxter algebra of weight $\theta \in A$. On
the underlying vector space we define the operation
\begin{equation}
\label{doubleprod}
    a *_R b := R(a)b + aR(b) + a \ \theta\ b,\qquad  a,b \in A.
\end{equation}
A simple calculation shows that for $[\theta,x]:=\theta x - x
\theta=0$, $x \in R(A)$ we find the result.

\begin{prop} \label{prop:doubleRB1}
Let $(A,R)$ be a Rota--Baxter algebra of weight $\theta \in A$.
Assume $[\theta,R(A)]=0$. The vector space underlying $A$ equipped
with the operation in (\ref{doubleprod}) turns out to be an
algebra. We denote this not necessarily unital algebra by $A_R$.
\end{prop}

\begin{cor} \label{cor:doubleRB2}
For $A$ being a Rota--Baxter algebra of scalar weight $\theta \in
\mathbb{K}$, the pair $(A_R,R$) is a Rota--Baxter algebra of
weight $\theta \in \mathbb{K}$.
\end{cor}

In general, for $(A,R)$ a Rota--Baxter algebra of weight $\theta
\in A$ we only have
\begin{equation}
    \label{eq:RBhom}
    R(a \ast_R b) = R(a) R(b).
\end{equation}
For $\tilde{R} = -\theta\id_A - R$, $\theta \in \mathbb{K}$
satisfying the Rota--Baxter relation (see
Proposition~\ref{prop:RBdual}) one shows that $\tilde{R}(a \ast_R
b) = - \tilde{R}(a) \tilde{R}(b)$ holds.

\begin{rmk}{\rm{In \cite{Baxter60} Baxter showed that
Spitzer's classical identity~\cite{Spitzer56} is satisfied for
commutative Rota--Baxter algebras $(A,R)$ of weight $\theta \in
A$, that is, for $a \in A$ fixed he proved in $A[[t]]$ that
$$
    \exp\Big(\sum_{n>0}\frac{r_n t^n}{n}\Big) = 1_A + \sum_{m > 0} a_m t^m
$$
where $r_n:=R((-\theta)^{n-1}a^n)$, and $a_1 = r_1 = R(a)$,
$$
    a_m = \sum_{(\lambda_1,\dots,\lambda_m)}
    \frac{r_1^{\lambda_1} \cdots r_m^{\lambda_m}}{1^{\lambda_1}2^{\lambda_2}
    \cdots m^{\lambda_m} \lambda_1! \cdots \lambda_m!}
    = \underbrace{R\bigl(R(R( \cdots R}_{m-times}(a)a)\dots a)a\bigr).
$$
The sum goes over all integer $m$-tuples
$(\lambda_1,\dots,\lambda_m)$, $\lambda_i \geq 0$ for which
$1\lambda_1+ \dots + m\lambda_m = m$.}}
\end{rmk}


\subsection{Weight $\theta=R(1_A)$ case: $TD$-algebra}
\label{ssect:TDalgebra}

Let $(A,R)$ be a Rota--Baxter algebra of weight $\theta \in A$.
Assume $\theta:=R(1_A)$. One verifies that
$$
    R(1_A)R(x)=R^2(x)=R(x)R(1_A)
$$
for $x \in A$, i.e., $R(1_A)$ lies in the center of $R(A)$ (recall
Proposition~\ref{prop:doubleRB1}). The element $R(1_A)$ plays a
special role in the classical theory of Rota--Baxter algebras. We
refer the reader to~\cite{NHB71,NHB73,NHB76,Mil66,Mil69}. As we
will see in the following sections the weight $\theta = - R(1_A)$
case distinguishes itself from the scalar weight case, $\theta \in
\mathbb{K}$.

\smallskip

The particular case of Rota--Baxter algebra of weight $\theta:=
-R(1_A)$ appeared in the work of one of us~\cite{Leroux04a} in the
context of Loday's dendriform algebras and Nijenhuis
algebras.

\smallskip

We will keep the name $TD$-algebra introduced in~\cite{Leroux04a}
referring to an algebra $A$ with a $\mathbb{K}$-linear map $P: A
\to A$ satisfying the $TD$-relation
\begin{equation}
    P(x)P(y) = P\bigl(P(x)y + xP(y)\bigr) - P(x \ P(1_A)\ y),\qquad \forall x,y \in A.
\label{eq:TD}
\end{equation}

For a commutative Rota--Baxter algebra of weight $\theta=-P(1)$ we
find that $\tilde{P}:=P(1_A)\id_A -P$ also is a Rota--Baxter map
of weight $\theta=-P(1)$. But we do not think of $\tilde{P}$ as a
$TD$-map since replacing $P$ by $\tilde{P}$ in relation
(\ref{eq:TD}) does not give a $TD$-relation, since
$\tilde{P}(1_A)=0$.

\begin{prop} \label{prop:RBtdNijenhuis}
Let $(A,P)$ be a $TD$-algebra. Then the map $\tilde{P}:=
P(1_A)\id_A - P$ is an associative Nijenhuis operator on $A$.
\end{prop}

\begin{proof}
\allowdisplaybreaks{
 \begin{eqnarray}
  \tilde{P}(x)\tilde{P}(y) +  \tilde{P}^2(xy) &=&
                             P(1_A)xP(1_A)y - P(1_A)xP(y) - P(x)P(1_A)y                               \nonumber\\
                           & & \hspace{5cm} + P(x)P(y) + \tilde{P}\bigl(P(1_A)xy - P(xy)\bigr)          \nonumber\\
                           &=& P\bigl(P(x)y+xP(y)-xP(1_A)y\bigr) +
                                                             P(1_A)\bigl(xP(1_A)y - xP(y) - P(x)y\bigr) \nonumber\\
                           & & \hspace{2.2cm} + P(1_A)P(1_A)xy - P(1_A)P(xy) - P(P(1_A)xy) + P^2(xy)  \nonumber\\
                           &=& \tilde{P}\bigl( xP(1_A)y - xP(y) - P(x)y + P(1_A)xy \bigr)               \nonumber\\
                           &=& \tilde{P}\bigl( \tilde{P}(x)y + x\tilde{P}(y) \bigr).
 \end{eqnarray}}
\noindent We used the relation $P(1_A)P(x)=P^2(x)=P(x)P(1_A)$,
i.e., $[P(1_A),P(x)]=0$, $x \in A$.
\end{proof}

Let $(A,P)$ be a $TD$-algebra. The vector space underlying $A$,
equipped with the operation
$$
    a *_P b = P(a)b + aP(b) - aP(1_A)b, \qquad a,b \in A,
$$
is associative~\cite{Leroux04b}. It is of $TD$-type only if
$P(1_A)=P(P(1_A))$. Hence, we do not find immediately the notion
of a double $TD$-algebra. But we make the following observation.

\begin{prop}
Let $(A,P)$ be a $TD$-algebra. The pair $(A_P,P)$ is an
associative Nijenhuis algebra.
\end{prop}

\begin{proof}
Let $a,b \in A$.
 \allowdisplaybreaks{
\begin{eqnarray*}
    P(a)*_P P(b) + P^2(a *_P b)
    &=& P(a)P^2(b) + P\bigl(P(a)P(b)\bigr) \\
    &=& P\bigl(P(a)P(b) + aP^2(b) - aP(1_A)P(b) \bigr) + P\bigl(P(a)P(b)\bigr) \\
    &=& 2P\bigl( P(a)P(b) \bigr)
     =  P(a)P^2(b) + P^2(a)P(b)\\
    &=& P\bigl(P(a) *_P b + a *_P P(b) \bigr).
\end{eqnarray*}}
\end{proof}


\subsection{Rota--Baxter type categories}
\label{ssect:RBcategories}

For the rest of this paper, only unital commutative Rota--Baxter
type algebras will be considered if not stated otherwise. In the
sequel we will refer to these Rota--Baxter type algebras as
$\textrm{{\bf{R}}}$-algebras in general.

The above results motivate the following point of view. We denote
by {\textbf{ComRBA}} the category of commutative unital
Rota--Baxter algebras of scalar weight $\theta \in \mathbb{K}$.
The category of commutative unital $TD$-algebras is denoted by
{\textbf{ComTDA}}. Finally, by \textbf{ComNA} we mean the category
of commutative Nijenhuis algebras.

Those categories have as objects unital commutative
$\textrm{{\bf{R}}}$-algebras and as morphisms so-called
$\textrm{{\bf{R}}}$-morphisms, defined as follows. Let $(A_1,X_1)$
and $(A_2,X_2)$ be two objects in
\textbf{Com}$\textrm{{\bf{R}}}$\textbf{A}. An
$\textrm{{\bf{R}}}$-algebra morphism ${\phi}:(A_1,X_1)
\longrightarrow (A_2,X_2)$ is a unital algebra homomorphism such
that the diagram
\begin{diagram}
    A_1        & \rTo^{\phantom{m} \phi \phantom{m}}       & A_2\\
   \dTo^{X_1}  &                                           & \dTo_{X_2}\\
    A_1        & \rTo^{\phantom{m} \phi \phantom{m}}       & A_2
\end{diagram}
\noindent commutes, i.e., ${\phi}X_1=X_2{\phi}$.

By \textbf{Inv}$\textrm{{\bf{R}}}$\textbf{A}, we mean the category
of involutive unital commutative $\textrm{{\bf{R}}}$-algebras with
all the involved morphisms being involutive ones.

Once defined, a natural question to ask is how these categories
are related and what are the respective free objects. This is the
purpose of the following sections. As we will see, the
construction of the free objects in the above categories of
$\textrm{{\bf{R}}}$-algebras relies on shuffle products and their
generalizations such as the quasi-shuffle product.

We will consider free $\textrm{{\bf{R}}}$-algebras over an algebra
$A$. We briefly recall the universal property for the free object
in the categories \textbf{Com}$\textrm{{\bf{R}}}$\textbf{A}. Let
$A$ be a unital commutative $\mathbb{K}$-algebra. The free unital
commutative $\textrm{{\bf{R}}}$-algebra over $A$ is by definition
a unital commutative $\textrm{{\bf{R}}}$-algebra denoted by
$(\textrm{{\bf{R}}}(A),P_A)$, equipped with a unital algebra
morphism $i_A: A \hookrightarrow \textrm{{\bf{R}}}(A)$, such that
for any $\textrm{{\bf{R}}}$-algebra $(B,X)$ together with a unital
algebra morphism $\phi: A \rightarrow B$, there exists a unique
$\textrm{{\bf{R}}}$-algebra morphism $\tilde{\phi}:
\textrm{{\bf{R}}}(A) \rightarrow B$ satisfying $\phi =
\tilde{\phi} \circ i_A$, i.e., the following diagram commutes
\begin{diagram}
    A          & \rInto^{\phantom{mm}i_A\phantom{mm}}   & \textrm{{\bf{R}}}(A)\\
               & \rdTo_{\phi\ }                         & \dTo_{{\tilde{\phi}} }\\
               &                                        & B
\end{diagram}

\noindent Baxter, Rota and Cartier gave several constructions of
free commutative Rota--Baxter algebras. The concept of free unital
commutative $\textrm{{\bf{R}}}$-algebra over an algebra $A$ was
used in \cite{GK00} to give another construction of the free
object in {\textbf{ComRBA}} in terms of the mixable shuffle
product, which is related to the quasi-shuffle
product~\cite{EGmixShuf}. Recently, the general construction of
free associative Rota--Baxter algebra was presented, using
generalized non-commutative shuffle-like products on decorated
planar rooted trees \cite{AguMor05,EFG05c}.

In \cite{KEF2} the free object in {\textbf{ComNA}} was explicitly
constructed by introducing another type of generalized shuffle
product, i.e. the right-shift shuffle. Below, we will briefly
present the construction of the free commutative $TD$-algebra over
an algebra $A$ using a similar generalization of the shuffle
product.\smallskip

The nature of the shuffle products (see next section), with the
mixable, or quasi-shuffle on the one side and the left- and
right-shift shuffle on the other side, will explain our separate
definition of {\textbf{ComRBA}} from {\textbf{ComTDA}} and
{\textbf{ComNA}}, respectively (c.f.
Subsection~\ref{ssect:shufflePord} and
paragraphs \ref{sssect:RightShift} -- \ref{sssect:LeftShift}).\\


\subsection{Rota--Baxter type algebras and dendriform algebras}
\label{ssect:dendlink}

In Section~\ref{sect:dendal} we introduced dendriform algebras. It
was Aguiar~\cite{Aguiar00} who first observed the close relation
between associative Rota--Baxter algebras, of weight zero, and
Loday's dendriform algebra. More general, if $R$ is a Rota--Baxter
operator of non-zero weight $\theta \in \mathbb{K}$ on an
associative algebra $A$, then the operations $x \prec y := x
R(y)$, $x \succ y:= R(x) y$, and $x \bullet y = \theta x y$
determine a tridendriform algebra structure on $A$. The
associative product in (\ref{def:double}) just reflects the double
Rota--Baxter product construction in
Subsection~\ref{sssect:doubleRBtype}. Recently, such connections
between associative scalar weight Rota--Baxter algebras and
dendriform algebras were further explored in
\cite{AguMor05,KEF1,EG04,EG05,Leroux04a}.


\section{Generalizations of the shuffle product}
\label{sect:shuffle}

For completeness we recall the definition of shuffle product and
introduce three generalizations. We first start with the non
recursive definitions of those products.

Let $V$ be a $\mathbb{K}$-vector space, the tensor algebra is
denoted by
$$
  T(V) := \bigoplus_{n\geq 0} V^{\otimes n},
$$
with $V^{\otimes 0}:=\mathbb{K}$. The augmentation ideal
$\bar{T}(V):=T(V)/\mathbb{K}1_\mathbb{K}=\bigoplus_{n \geq 1}
V^{\otimes n}$. Recall that the vector space $T(V)$ is naturally
equipped with a grading, $l(u)=|u|=n$ for $u \in V^{\otimes n}$.\\


\subsection{Shuffle product}
\label{ssect:shufflePord}

The shuffle product on $T(V)$ starts with the shuffles of
permutations~\cite{Reutenauer93,Sweedler69}. For $m,n \in
\mathbb{N}_+$, define the set of {\bf{$(m,n)$-shuffles}} by
\[
  S(m,n)= \left \{ \sigma\in S_{m+n}
          \begin{array}{ll}
                {} \\ {}
          \end{array} \right .
\left |
\begin{array}{l}
 \sigma^{-1}(1)<\sigma^{-1}(2)<\ldots<\sigma^{-1}(m),\\
 \sigma^{-1}(m+1)<\sigma^{-1}(m+2)<\ldots<\sigma^{-1}(m+n)
\end{array}
\right \}.
\]
Here $S_{m+n}$ is the symmetric group on $m+n$ letters.

For $a = a_1 \otimes \ldots \otimes a_m \in V^{\otimes m}$, $b =
b_1 \otimes \ldots \otimes b_n \in V^{\otimes n}$ and $\sigma \in
S(m,n)$, the element
\[
  \sigma (a \otimes b) := u_{\sigma(1)} \otimes u_{\sigma(2)} \otimes
                                            \ldots \otimes u_{\sigma(m+n)} \in V^{\otimes (m+n)},
\]
where
\[
    u_k = \left \{ \begin{array}{ll}
                     a_k,     & \quad   1 \leq k \leq m,\\
                     b_{k-m}, & \quad m+1 \leq k \leq m+n,
                   \end{array}
    \right.
\]
is called a {\bf{shuffle}} of $a$ and $b$.

The {\bf{shuffle product}} $(a_1 \otimes \ldots \otimes a_m)
\:\sha\: (b_1 \otimes \ldots \otimes b_n)$ is the sum over all
shuffles of $a_1 \otimes \ldots \otimes a_m$ and $b_1 \otimes
\ldots \otimes b_n$
\begin{equation}
  a \:\sha\: b := \sum_{\sigma \in S(m,n)} \sigma(a \otimes b).
 \label{eq:shuffle}
\end{equation}

\begin{example}{\rm{
    $a_1 \:\sha\: (b_1 \otimes b_2) = a_1 \otimes b_1 \otimes b_2
                                      + b_1 \otimes a_1 \otimes b_2 + b_1 \otimes b_2 \otimes a_1.$}}
\end{example}

Also, by convention, $a \:\sha\: b$ is the scalar product if
either $m=0$ or $n=0$, that is, if $a$ or $b$ is in $V^{\otimes 0}
=\mathbb{K}$.

\begin{prop} \label{prop:shuffleAlg}
The bilinear composition $\:\sha\:$ defines an associative,
commutative bilinear product on $T(V)$, making
$Sh(V):=(T(V),\sha,1_{\mathbb{K}})$ into a unitary
$\mathbb{K}$-algebra with unit $1_\mathbb{K} \in V^{\otimes 0}$.
\end{prop}

Moreover, one verifies the following proposition.

\begin{prop} \label{lem:shufRB}
For every $u \in V$, the $\mathbb{K}$-linear map $P^{(u)}_V: Sh(V)
\to Sh(V)$, $P^{(u)}_V(a):= u \otimes a$, $a \in T(V)$, is a
Rota--Baxter operator of weight zero.
\end{prop}

\begin{proof}
The proof reduces to a simple calculation, using the recursive
definition of the shuffle product given further below.
\end{proof}


\subsection{Generalized shuffle products}
\label{ssect:shufflePord}

We now work with the tensor algebra, $T(A)$, over an associative
unital $\mathbb{K}$-algebra $A$ and equip it with three
generalized shuffle products. We give explicit formulae for these
products before stating their recursive definitions. Here, we
follow the terminology in \cite{EGmixShuf,GK00}.

Let us start with a $(m,n)$-shuffle $\sigma \in S(m,n) \subset
S_{m+n}$. A pair of indices $(k,k+1)$, $1 \leq k < m+n$, is called
an {\bf{admissible pair for}} $\sigma$ if $\sigma(k) \leq m <
\sigma(k+1)$. Denote $\mathcal{T}^\sigma$ for the set of all
admissible pairs for the $(m,n)$-shuffle $\sigma$. For a subset
$T$ of $\mathcal{T}^\sigma$, call the pair $(\sigma,T)$ a
{\bf{mixable $(m,n)$-shuffle}}. Let $|T|$ be the cardinality of
$T$. By convention, $(\sigma,T) = \sigma$ if $T=\emptyset$. Denote
\[
     {\bar{S}} (m,n)=\{ (\sigma,T) \big{|} \sigma \in S(m,n), T \subset \mathcal{T}^\sigma\}
\]
for the set of mixable $(m,n)$-shuffles.

Let $(A,m_A,1_A)$ be a unital associative $\mathbb{K}$-algebra.
For $n \geq 0$ let $A^{\otimes n}$ be the $n$-th tensor power of
$A$ over $\mathbb{K}$ with the convention that $A^{\otimes
0}=\mathbb{K}$ and define the tensor module
$$
  T(A) =\bigoplus_{n\geq 0} A^{\otimes n}.
$$
For two tensors $a = a_1 \otimes \ldots \otimes a_m \in
A^{\otimes m}$, $b = b_1 \otimes \ldots \otimes b_n \in A^{\otimes
n}$ and $(\sigma,T)\in {\bar{S}}(m,n)$, $T=\emptyset$ the element
\[
   \sigma(a \otimes b;\emptyset) := \sigma (a \otimes b) = u_{\sigma(1)} \otimes u_{\sigma(2)} \otimes
                                                           \ldots \otimes u_{\sigma(m+n)} \in A^{\otimes (m+n)},
\]
where
\[
    u_k:=\begin{cases}        a_k,     & \quad   1 \leq k \leq m,\\
                          b_{k-m},     & \quad m+1 \leq k \leq m+n,
          \end{cases}
\]
was called a shuffle of $a$ and $b$. Now, we extend this notion as
follows. The element
\[
  \sigma_q(a \otimes b; T) := u_{\sigma(1)}{\hat{\otimes}_q} u_{\sigma(2)}
                                                   {\hat{\otimes}_q} \ldots {\hat{\otimes}_q} u_{\sigma(m+n)},
\]
where for each pair $(k,k+1)$, $1 \leq k < m+n$,
\[
    u_{\sigma(k)} {\hat{\otimes}_q} u_{\sigma(k+1)}
                        :=\left\{\begin{array}{ll}
                         \theta m_A \bigl(u_{\sigma(k)} \otimes u_{\sigma(k+1)}\bigr) \in A,
                                                        & \quad (k,k+1)\in T\ \mathrm{and}\ q=\theta\\
                         -1_A \otimes m_A \bigl(u_{\sigma(k)} \otimes u_{\sigma(k+1)}\bigr) \in A \otimes A,
                                                        & \quad (k,k+1)\in T\ \mathrm{and}\ q=\mathrm{r}\\
                         -m_A \bigl(u_{\sigma(k)} \otimes u_{\sigma(k+1)}\bigr) \otimes 1_A \in A \otimes A,
                                                        & \quad (k,k+1)\in T\ \mathrm{and}\ q=\ell\\
                                  u_{\sigma(k)} \otimes u_{\sigma(k+1)} \in A \otimes A,
                                                        & \quad (k,k+1) \not \in T,
                           \end{array} \right.
\]
is called {\bf{mixable shuffle}} ($q=\theta$) of scalar weight
$\theta \in \mathbb{K}$, or {\bf{right-shift}} respectively
{\bf{left-shift}} shuffle ($q=\mathrm{r}$ resp. $q=\ell$) of the
tensors $a$ and $b$.

Define, for $a \in A^{\otimes m}$ and $b \in A^{\otimes n}$ as
above and $q \in\{\theta,\mathrm{r},\ell\}$
\begin{equation}
     a \bullet^q b\:= \sum_{(\sigma,T)\in {\bar{S}} (m,n)} \sigma_q(a\otimes b;T)
      \in
      \left\{
       \begin{array}{ll}
        \bigoplus_{k \leq m+n} A^{\otimes k}, & q=\theta,\\
        A^{\otimes m+n},                      & q=\mathrm{r}, \ell,
       \end{array} \right.
\label{eq:geneShuffle}
\end{equation}
to be the quasi- or mixable shuffle product of weight $\theta$,
and the right-shift shuffle and left-shift shuffle product,
respectively. As usual, for $a$ or $b$ being scalars, we define $a
\bullet^q b := ab$.

We see immediately that the case $q=\theta =0$ reduces to the
ordinary shuffle product, $\bullet^0 = \sha$.
\begin{example}{\rm{ A simple example might help to keep track
of the products.
$$
    a_1 \left\{\begin{array}{c} \bullet^\theta \\
                                \bullet^\mathrm{r}\\
                                \bullet^\ell
               \end{array}\right\} (b_1 \otimes b_2)
               = a_1 \:\sha\: (b_1 \otimes b_2) +
               \left\{\begin{array}{c}
                              + \theta m_A(a_1 \otimes b_1) \otimes b_2 + \theta b_1 \otimes m_A(a_1 \otimes b_2)            \\
                              - 1_A \otimes m_A(a_1  \otimes b_1) \otimes b_2 - b_1 \otimes 1_A \otimes m_A(a_1 \otimes b_2)  \\
                              - m_A(a_1 \otimes b_1) \otimes 1_A  \otimes b_2 - b_1 \otimes m_A(a_1 \otimes b_2) \otimes 1_A
               \end{array}\right\}
$$}}
\end{example}

Observe that all three products have in common that they consist
of two parts, i.e.,
$$
     a \bullet^q b = a \:\sha\: b + M^q(a,b),
$$
where $a \:\sha\: b$ corresponds to the proper $(m,n)$-shuffles
$\sigma_q(a\otimes b;\emptyset) \in \bar{S}(m,n)$ and $M^q(a,b)$
is a sum of shuffles modified by contractions of adjacent pairs of
letters in $A$. Hence, the product is not grade (i.e. length)
preserving in the case $q=\theta$. However, the two other cases,
$q=\ell,\mathrm{r}$, preserve the grading by length due to the
extra left- or right-{\it{shift}}, respectively, by an element of
$A$. More important is the fact that --associativity and--
commutativity of the above shuffle products only depends on the
product in $A$ and its properties, i.e. whether $A$ is commutative
or not.

As in the case of the shuffle product, it has been shown for the
cases $q=\theta$ \cite{GK00} and $q=\mathrm{r}$ \cite{KEF2} that
the operation $\bullet^q$ extends to an associative binary
operation on $T(A)$ making $\sha^{q}(A):=(T(A),\bullet^q )$ into
an associative unital algebra with unit $1_\mathbb{K}$.
Commutativity of $A$ implies commutativity of these generalized
shuffle products.\\


\subsection{Recursive definition of generalized shuffle products and the shift-maps}
\label{ssect:recursiveShuffle}

In the case of the classical shuffle, the mixable, or
quasi-shuffle, and the right-shift shuffle product recursive
formulae appeared already in several places,
e.g.~\cite{BowmanBrad02,KEF2,EGmixShuf,GK00,Hoffman00}.
Nevertheless, we will briefly recall them for the sake of
completeness.\\


\subsubsection{Recursive shuffle product}
\label{sssect:recShuffle}

Let us begin with the classical shuffle product on the tensor
module $T(V)$, $V$ being a $\mathbb{K}$-vector space. We choose
two elements $a_1 \otimes \cdots \otimes a_m \in V^{\otimes m}$
and $b_1 \otimes \cdots \otimes b_n \in V^{\otimes n}$, and define
 \allowdisplaybreaks{
\begin{eqnarray*}
      a_0 1_\mathbb{K} \:\sha\: (b_1\otimes b_2 \otimes \ldots \otimes b_n)  &
                        =&  a_0 b_1 \otimes b_2\otimes \ldots \otimes b_n, \\
      (a_1 \otimes a_2 \otimes \ldots \otimes a_m) \:\sha\: b_0 1_\mathbb{K} &
                        =&  b_0 a_1 \otimes a_2 \otimes \ldots \otimes a_m,
\end{eqnarray*}}
for scalars $a_0,b_0$ in $V^{\otimes 0}={\mathbb{K}}$, and in
general
 \allowdisplaybreaks{
\begin{eqnarray}
      \lefteqn{(a_1\otimes  \ldots \otimes a_m) \:\sha\: (b_1 \otimes \ldots \otimes b_n) } \notag \\
             &:=& a_1 \otimes \big ((a_2\otimes \ldots \otimes a_m) \:\sha\:
                                          (b_1\otimes \ldots \otimes b_n)\big ) \notag\\
             & & \;\;\qquad\qquad \ +\ b_1 \otimes \big ((a_1 \otimes \ldots \otimes a_m)
                        \:\sha\: (b_2\otimes \ldots \otimes b_n)\bigr),\qquad a_i,\ b_j \in V.
      \label{eq:rshuf1}
\end{eqnarray}}
With the recursive definition of the shuffle product the proof of
Proposition~\ref{lem:shufRB} follows by showing that
$P^{(u)}_V(W)=P^{(u)}(W) := u \otimes W$ on $Sh(V)$ defines a
weight zero Rota--Baxter map,
\allowdisplaybreaks{
\begin{eqnarray*}
    P^{(u)}(U) \:\sha\: P^{(u)}(W) = (u \otimes U) \:\sha\: (u \otimes W)
                                   &=& u \otimes (U \:\sha\: u \otimes W) + u \otimes (u \otimes U \:\sha\: W)\\
                                   &=& P^{(u)} \bigl(P^{(u)}(U) \:\sha\: W + U \:\sha\: P^{(u)}(W)\bigr).
\end{eqnarray*}}


\subsubsection{Recursive quasi- or mixable shuffle product}
\label{sssect:quasishuffle}

Now, replace the vector space $V$ in $T(V)$ by an associative but
not necessarily commutative $\mathbb{K}$-algebra $A$ with unit
$1_A$. We denote by $[a;b]$ the associative algebra product for
$a,b \in A$. To streamline our presentation we introduce the
following notation. Tensor strings in $T(A)$ are denoted by
capital letters, whereas lower case letters denote elements in
$A$. Recall, that $1_\mathbb{K}$ denotes the empty word in $T(A)$. \\

Let $\theta \in \mathbb{K}$ be a scalar parameter. Then the
{\bf{quasi-}} or {\bf{mixable shuffle}} product of weight $\theta$
on $T(A)$ has the recursive definition~\cite{GK00,Hoffman00}
\begin{equation}
  \label{eq:mixShuf}
  a \otimes U \ \bullet^\theta \ b \otimes W := a \otimes (U \bullet^\theta b \otimes W)
                                                 + b \otimes (a \otimes U \bullet^\theta W)
                                                  + \theta [a;b] \otimes (U \bullet^\theta W),
\end{equation}
and $U \bullet^\theta k1_\mathbb{K} := kU =:k1_\mathbb{K}
\bullet^\theta U$, for $k \in \mathbb{K}$. Recall that
$\sha^\theta(A):=(T(A),\bullet^\theta,1_\mathbb{K})$ is an
associative --commutative-- unital algebra. For $\theta=0$ we are
back to the shuffle algebra.\\


\subsubsection{Recursive right-shift shuffle product}
\label{sssect:RightShift}

The recursive formula for the {\bf{right-shift shuffle}} product
on $T(A)$ is simply given by
\begin{equation}
\label{eq:rightShuf}
  a \otimes U \bullet^\mathrm{r} b \otimes V = a \otimes (U \bullet^\mathrm{r} b \otimes V)
                                                 + b \otimes (a \otimes U \bullet^\mathrm{r} V)
                                                    - 1_A \otimes [a;b] \otimes (U \bullet^\mathrm{r} V),
\end{equation}
$U \bullet^\mathrm{r} k1_\mathbb{K} := kU =: k1_\mathbb{K}
\bullet^\mathrm{r} U$, for $k \in \mathbb{K}$. And
$\sha^\mathrm{r}(A):=(T(A),\bullet^\mathrm{r},1_\mathbb{K})$ is an
associative --commutative-- unital algebra~\cite{KEF2}.\\


\subsubsection{Recursive left-shift shuffle product}
\label{sssect:LeftShift}

Finally, we introduce a recursive formula for the {\bf{left-shift
shuffle}} product on $T(A)$.
\begin{equation}
\label{eq:leftShuf}
   a \otimes U \bullet^\ell b \otimes V = a \otimes (U \bullet^\ell b \otimes V)
                                            + b \otimes (a \otimes U \bullet^\ell V)
                                               - [a;b] \otimes 1_A \otimes (U \bullet^\ell V),
\end{equation}
$U \bullet^\ell k 1_\mathbb{K}:= kU =: k 1_\mathbb{K} \bullet^\ell
U$, for any $k \in \mathbb{K}$. And
$\sha^\ell(A):=(T(A),\bullet^\ell,1_\mathbb{K})$ is an associative
--commutative-- unital algebra. We postpone the proof of this
statement to the next section.\\


\subsection{Rota--Baxter type algebras and generalized shuffle products}
\label{ssect:RotaBaxterTypeAlg}

We introduce two linear maps on $T(A)$ in terms of the algebra
unit $1_A$. First, we define the {\bf{right-shift map}}
$$
    P_A(U) := 1_A \otimes U, \qquad U \in T(A)/\mathbb{K}1_{\mathbb{K}},
$$
and the {\bf{left-shift map}}
$$
    Q_A(U) := U \otimes 1_A ,\qquad U \in T(A)/\mathbb{K}1_{\mathbb{K}}.
$$
Both maps commute and we extend them to $T(A)$ by defining
$P_A(1_\mathbb{K}) := 1_A =: Q_A(1_{\mathbb{K}})$. These two maps
are crucial with respect to the construction of free objects in
the categories \textbf{Com}$\mathrm{\bf{R}}$\textbf{A}
of commutative unital $\mathrm{{\bf{R}}}$-algebras.\\

Let $T(A)$ be equipped with one of the above shuffle products
$\bullet^q$, $q \in\{\theta, \mathrm{r},\ell\}$. We regard
$$
    T^+(A) := A \otimes T(A),
$$
with the standard product
\begin{equation}
\label{extShuffle}
    a \otimes U \:\bar{\bullet}^q\: b \otimes V := [a;b] \otimes (U \:\bullet^{q}\: V).
\end{equation}

First, observe that $A$ naturally imbeds into $T^{+}(A)$ by the
map $i_A : A \hookrightarrow T^{+}(A)$
\begin{equation}
\label{inclusion}
    a \mapsto i_A(a) := a \otimes 1_{\mathbb{K}}.
\end{equation}
Secondly, this inclusion is an algebra morphism
 \allowdisplaybreaks{
\begin{eqnarray*}
    i_A\bigl([a; b]\bigr) &=&  [a ; b] \otimes 1_{\mathbb{K}}
                           =\bigl( a \otimes 1_{\mathbb{K}}  \bigr)
                                             \:\bar{\bullet}^q\:
                                   \bigl( b \otimes 1_{\mathbb{K}}\bigr) = i_A(a) \:\bar{\bullet}^q\: i_A(b).
\end{eqnarray*}}
One shows that $i_A(1_A) = 1_A \otimes 1_{\mathbb{K}}$ is the unit
for $\bar{\bullet}^q$, $q \in \{\theta, \mathrm{r},\ell\}$. We
extend the right- and left-shift map, $P_A$ respectively $Q_A$, to
$T^{+}(A)$ by defining $P_A(a \otimes 1_{\mathbb{K}}):= 1_A
\otimes a$, analogously for $Q_A$.

With the exception of Section~\ref{sect:Rbialg}, we will omit the
tensor product by $1_{\mathbb{K}}$ in the notation for elements in
$i_A(A) \subset T^{+}(A)$ in the sequel.\\


\subsubsection{The free Rota--Baxter and Nijenhuis algebra}
\label{sssect:freeRBandNA}

The algebra $\sha^\theta(A)=(T(A),\bullet^\theta,1_\mathbb{K})$
together with the map $P_A$ defines a unital Rota--Baxter algebra
of scalar weight $\theta$. This follows from the recursive
definition of $\bullet^\theta$
 \allowdisplaybreaks{
\begin{eqnarray*}
    P_A(U) \bullet^\theta P_A(V) &=& 1_A \otimes U \bullet^\theta 1_A \otimes V \\
                                 &=& 1_A \otimes \bigl(U \bullet^\theta 1_A \otimes V \bigr)
                                      + 1_A \otimes \bigl(1_A \otimes U \bullet^\theta V \bigr)
                                       + \theta [1_A;1_A] \otimes (U \bullet^\theta V )\\
                                 &=& P_A \bigl(U \bullet^\theta P_A(V)\bigr)
                                      + P_A \bigl(P_A(U) \bullet^\theta V\bigr)
                                       + \theta  P_A(U \bullet^\theta V).
\end{eqnarray*}}

Moreover, if $A$ is commutative, then the free unital commutative
Rota--Baxter algebra of scalar weight $\theta$ over $A$ is simply
$T^+(A):=A \otimes T(A)$ with corresponding product
(\ref{extShuffle}) and unit $1_A$. We denote it by
$\overline{\sha}^\theta(A):=(T^+(A),
\bar{\bullet}^\theta,P_A)$~\cite{GK00}.\medskip

The map $P_A$ in connection with the right-shift shuffle product
gives birth to an associative Nijenhuis algebra $(T(A),
\bullet^{\mathrm{r}}, P_A)$. Again, this is readily seen by the
recursive definition of the right-shift shuffle product
 \allowdisplaybreaks{
\begin{eqnarray*}
    P_A(U) \bullet^\mathrm{r} P_A(V) &=& 1_A \otimes U \bullet^\mathrm{r} 1_A \otimes V \\
                                     &=& 1_A \otimes \bigl(U \bullet^\mathrm{r} 1_A \otimes V \bigr)
                                           + 1_A \otimes \bigl(1_A \otimes U \bullet^\mathrm{r} V \bigr)
                                            - 1_A \otimes [1_A;1_A] \otimes (U \bullet^\mathrm{r} V )\\
                                     &=& P_A \bigl(U \bullet^\mathrm{r} P_A(V)\bigr)
                                           + P_A \bigl(P_A(U) \bullet^\mathrm{r} V\bigr)
                                            -  P^2_A (U \bullet^\mathrm{r} V).
\end{eqnarray*}}

For $A$ a commutative algebra, $T^+(A):= A \otimes T(A)$ with unit
$1_A$ and product $\bar{\bullet}^\mathrm{r}$~(\ref{extShuffle})
together with the map $P_A$ is the free object, denoted by
$\overline{\sha}^\mathrm{r}(A):=(T^+(A),
\bar{\bullet}^\mathrm{r},P_A)$, in the category {\textbf{ComNA}}~\cite{KEF2}.\\


\subsubsection{The free Nijenhuis algebra is a $TD$-algebra}
\label{sssect:NijenhuisTD}

We now show that the free unital commutative Nijenhuis algebra
$\overline{\sha}^\mathrm{r}(A)$ over $A$ contains a $TD$-structure
in terms of the left-shift map $Q_A$.

\begin{lem}
Let $A$ be a unital associative algebra. For the right-shift
shuffle $\bullet^\mathrm{r}$ on $T(A)$, the equalities
$$
    X \bullet^\mathrm{r} 1_A^{\otimes n} =  X \otimes 1_A^{\otimes n}
                                         =  1_A^{\otimes n} \bullet^\mathrm{r} X
$$
for $n \in \mathbb{N} $ hold.
\end{lem}

\noindent Observe that we do not assume commutativity for $A$. The
proof follows by induction. This observation implies that the
subalgebra spanned by the set $I:=\{1_A^{\otimes n} \ | \ n\in
\mathbb{N} \} \subset Z(\sha^\mathrm{r}(A))$, where
$Z(\sha^\mathrm{r}(A))$ denotes the center of
$\sha^\mathrm{r}(A)$, which is generated by the center of $A$.

\begin{prop}
Consider $\overline{\sha}^\mathrm{r}(A)$ over a commutative
algebra $A$ with unit $1_A$. The map $Q_A$ is a $TD$-operator on
$\overline{\sha}^\mathrm{r}(A)$.
\end{prop}

\begin{proof}
We must show that $Q_A$ satisfies
$$
    Q_A(a\otimes X) \bar{\bullet}^\mathrm{r} Q_A(b\otimes Y) +
    Q_A \bigl(a\otimes X \bar{\bullet}^\mathrm{r} Q_A(1_A) \bar{\bullet}^\mathrm{r} b\otimes Y\bigr)
    = Q_A \bigl(Q_A(a\otimes X) \bar{\bullet}^\mathrm{r} b\otimes Y + a\otimes X \bar{\bullet}^\mathrm{r} Q_A(b\otimes Y)\bigr).
$$
But this follows from the last lemma, which completes the
proof.
\end{proof}

\noindent Therefore, the free  unital commutative Nijenhuis
algebra $\overline{\sha}^\mathrm{r}(A)$ is also a $TD$-algebra.
Moreover it has the structure of a commutative
dendriform-Nijenhuis algebra, see~\cite{Leroux04b}.

However, observe that the commutativity for $A$ can be dropped. If
this is the case, $(T(A)^+,\bar{\bullet}^\mathrm{r},P_A)$ is still
a Nijenhuis algebra as well as a $TD$-algebra via the map $Q_A$,
but one cannot think of it as the free object. In the following
section, we construct explicitly the free unital commutative
$TD$-algebra over $A$.\\


\subsubsection{The free commutative unital $TD$-algebra}
\label{sssect:freeTD}

\noindent To construct the free object in \textbf{ComTDA}, we
proceed as follows.

\begin{prop} \label{prop:leftshuffleAlg}
Let $A$ be an associative algebra (not necessarily commutative)
with unit $1_A$. Denote by $[-;-]$ its product. On $T(A)$ with the
left-shift shuffle $\bullet^\ell$ (\ref{eq:leftShuf}) we find
\begin{enumerate}
    \item For any $n \in \mathbb{N}$ and $X \in T(A)$ we have
\begin{equation}
    \label{eq:lID}
     X \bullet^\ell 1_A^{\otimes n} =  1_A^{\otimes n} \otimes X
                                    =  1_A^{\otimes n} \bullet^\ell X.
\end{equation}

    \item The left-shift shuffle is associative with the empty word $1_\mathbb{K}$ as unit, and
    commutative if and only if the algebra $A$ is commutative.
\end{enumerate}
\end{prop}

\begin{proof}
Recall that $a \bullet^\ell b := a \otimes b + b \otimes a - [a;b]
\otimes 1_A$ for $a,b \in A$, hence commutativity follows if and
only if $[-;-]$ on $A$ is commutative. The rest follows by
induction.
\end{proof}

\begin{prop} \label{prop:TDrelShuffle}
The right-shift operator $P_A$ on $(T(A),\bullet^\ell)$ fulfils
the $TD$-relation. That is,
 \allowdisplaybreaks{
\begin{eqnarray*}
    P_A(a \otimes X) \bullet^\ell P_A(b \otimes Y) &=& P_A \bigl(a \otimes X \bullet^\ell P_A(b \otimes Y)   \bigr)
                                                          +  P_A \bigl( P_A(a \otimes X) \bullet^\ell b \otimes Y  \bigr)      \\
                      & & \hspace{4cm} - P_A \bigl(a \otimes X  \bullet^\ell P_A(1_\mathbb{K}) \bullet^\ell b \otimes Y \bigr).
\end{eqnarray*}}
Hence, $\sha^\ell(A):=(T(A),\bullet^\ell,P_A)$ is a unital
$TD$-algebra.
\end{prop}

\begin{proof}
Remember that $P_A(1_\mathbb{K})=1_A$.
 \allowdisplaybreaks{
\begin{eqnarray}
 \lefteqn{ P_A(a \otimes X) \bullet^\ell P_A(b \otimes Y) = 1_A \otimes a \otimes X \bullet^\ell 1_A \otimes b \otimes Y} \nonumber\\
                  &=& 1_A \otimes (a \otimes X \bullet^\ell 1_A \otimes b \otimes Y)
                             + 1_A \otimes (1_A \otimes a \otimes X \bullet^\ell b \otimes Y)
                                                - [1_A; 1_A] \otimes 1_A \otimes ( a \otimes X \bullet^\ell b \otimes Y) \nonumber \\
                  &=& P_A\bigl(a \otimes Y \bullet^\ell P_A(b \otimes Y)\bigr)
                             + P_A \bigl( P_A(a \otimes X) \bullet^\ell b \otimes Y \bigr)
                              - 1_A \otimes 1_A \otimes ( a \otimes X\bullet^\ell  b \otimes Y),
    \label{eq:TDNijenhuis}
\end{eqnarray}}
which completes the proof using relation~(\ref{eq:lID}).
\end{proof}

\begin{rmk}{\rm{Proposition~\ref{prop:leftshuffleAlg} immediately
implies that $\sha^\ell(A)=(T(A),\bullet^\ell,P_A)$ is also a
Nijenhuis algebra. Indeed, using associativity of $\bullet^\ell$
and the first item of Proposition~\ref{prop:leftshuffleAlg}, we
get,
 \allowdisplaybreaks{
\begin{eqnarray*}
    P_A(a \otimes X \bullet^\ell P_A(1_\mathbb{K}) \bullet^\ell b \otimes Y) &=&
                        1_A \otimes \bigl((a \otimes X \bullet^\ell 1_A) \bullet^\ell b \otimes Y \bigr) \\
        &=& 1_A \otimes \bigl( (1_A \bullet^\ell  a \otimes X) \bullet^\ell b \otimes Y \bigr) \\
        &=& 1_A \otimes \bigl(  1_A \bullet^\ell (a \otimes X \bullet^\ell b \otimes Y) \bigr) \\
        &=& 1_A \otimes 1_A \otimes (a \otimes X \bullet^\ell b \otimes Y)
         =  P_A^2(a \otimes X \bullet^\ell b \otimes Y).
\end{eqnarray*}}
Observe that commutativity of $A$, and hence of $\sha^\ell(A)$, is
not needed here.}}
\end{rmk}

\noindent However, it is not true a priori that the Nijenhuis
algebra $\sha^\mathrm{r}(A)$ is also of $TD$-type. Moreover,
observe that $P_A(a \otimes X) \bullet^\ell P_A(b \otimes Y) =
P_A^2(a \otimes X \bullet^\ell b \otimes Y)$ and that,
$$
    P_A(a \otimes X) \bullet^\ell P_A(b \otimes Y) = P_A \bigl(a \otimes X \bullet^\ell P_A(b \otimes Y) \bigr)
                                                   = P_A \bigl( P_A(a \otimes X) \bullet^\ell b \otimes Y \bigr),
$$
showing that the right-shift operator $P_A$ in fact is an average
operator on $(T(A),\bullet^\ell)$, see
\cite{Aguiar00,Leroux03,EG05}.

\begin{prop}
Let $A$ be an associative algebra with unit $1_A$. The associative
$\mathbb{K}$-algebra $(T^+(A), \bar{\bullet}^{\ell},P_A)$ is a
$TD$-algebra, with unit $1_A$.
\end{prop}

\begin{proof}
Recall that the product $\bar{\bullet}^\ell$ on $T^+(A)$ is given by
$$
  (a \otimes X) \:\bar{\bullet}^{\ell}\: (b \otimes Y) = [a;b] \otimes (X \bullet^\ell Y),
$$
and that it is commutative and associative for $A$ being so. We
find that
 \allowdisplaybreaks{
\begin{eqnarray*}
 P_A(a \otimes X) \:\bar{\bullet}^{\ell}\: P_A(b \otimes Y)
           &=& 1_A \otimes \bigl( a \otimes X \bullet^{\ell} b \otimes Y \bigr)        \\
           &=& 1_A \otimes \bigl( a \otimes (X \bullet^{\ell} b \otimes Y)
                                + b \otimes (a \otimes X \bullet^{\ell} Y)
                                 - [a;b] \otimes 1_A \otimes (X \bullet^{\ell} Y) \bigr),
\end{eqnarray*}}
and
 \allowdisplaybreaks{
\begin{eqnarray*}
    P_A \bigl( P_A(a \otimes X) \:\bar{\bullet}^{\ell}\: b \otimes Y \bigr)
                            &=& 1_A \otimes \bigl( b \otimes (a \otimes X \bullet^{\ell} Y) \bigr),  \\
    P_A \bigl( a \otimes X  \:\bar{\bullet}^{\ell}\: P_A(b \otimes Y) \bigr)
                            &=& 1_A \otimes \bigl(a \otimes (X \bullet^{\ell} b \otimes Y) \bigr),  \\
   -P_A \bigl((a \otimes X) \:\bar{\bullet}^{\ell}\: P_A(1_A) \:\bar{\bullet}^{\ell}\: (b \otimes Y) \bigr)
                            &=& - 1_A \otimes \bigl( a \otimes (X \bullet^{\ell} 1_A)  \:\bar{\bullet}^{\ell}\: (b \otimes Y) \bigr)\\
                            &=& - 1_A \otimes [a;b] \otimes (X \:\bullet^{\ell}\: 1_A \:\bullet^{\ell}\: Y),\\
                            &=& - 1_A \otimes [a;b] \otimes (1_A \:\bullet^{\ell}\: X \:\bullet^{\ell}\: Y),\\
                            &=& - 1_A \otimes [a;b] \otimes  1_A \otimes (X \:\bullet^{\ell}\: Y).
\end{eqnarray*}}
Therefore,
$$
    P_A(a \otimes X) \bar{\bullet}^{\ell} P_A(b \otimes Y) = P_A \bigl( a \otimes X \bar{\bullet}^{\ell} P_A(b \otimes Y)
                                            +  P_A(a \otimes X) \bar{\bullet}^{\ell} b \otimes Y
                                            - a \otimes X  \bar{\bullet}^{\ell} P_A(1_A) \bar{\bullet}^{\ell} b \otimes Y \bigr).
$$
\end{proof}

\begin{thm} \label{thm:freeTD}
Let $A$ be a unital commutative algebra with unit $1_A$. The free
unital commutative $TD$-algebra over $A$ is given by the
$TD$-algebra $\overline{\sha}^\ell(A) := (T^+(A),
\bar{\bullet}^{\ell}, P_A)$.
\end{thm}

\begin{proof}
The previous Proposition shows for $A$ being commutative that
$(T^+(A), \bar{\bullet}^{\ell}, P_A)$ is a commutative
$TD$-algebra with unit $1_A$. Recall the canonical inclusion $i_A:
A \hookrightarrow \overline{\sha}^\ell(A)$~(\ref{inclusion}). Let
$(T,P)$ be a $TD$-algebra and $\phi: A \rightarrow T$ be an
algebra morphism. Define $\tilde{\phi}: T^+(A) \rightarrow T$ by
$\tilde{\phi} = \phi$ on $A$ and by induction as follows,
$$
 \tilde{\phi}(a_1 \otimes a_2 \otimes \ldots \otimes a_n) := \phi(a_1) P\bigl(\tilde{\phi}(a_2 \otimes \ldots \otimes a_n)\bigr).
$$
We have to prove that $\tilde{\phi}$ is a $TD$-algebra morphism.
Observe first that,
$$
    a_1 \otimes a_2 \otimes \ldots \otimes a_n = a_1 \bar{\bullet}^{\ell} P_A(a_2 \otimes \ldots \otimes a_n),
$$
leading by induction to $a_1 \otimes a_2 \otimes \ldots \otimes
a_n = a_1 \bar{\bullet}^{\ell} P_A\bigl(a_2 \bar{\bullet}^{\ell}
P_A( \ldots P_A(a_{n-1} \bar{\bullet}^{\ell}
P_A(a_n))\ldots)\bigr).$ Therefore,
 \allowdisplaybreaks{
\begin{eqnarray*}
    \tilde{\phi}\bigl(P_A( a_1 \otimes a_2 \otimes \ldots \otimes a_n) \bigr)
                        &=&  \tilde{\phi}\bigl(1_A \otimes  a_1 \otimes a_2 \otimes \ldots \otimes a_n \bigr), \\
                        &:=& \phi(1_A) P\bigl( \tilde{\phi}(a_1 \otimes a_2 \otimes \ldots \otimes a_n) \bigr),\\
                        &=&  P\bigl( \tilde{\phi}(a_1 \otimes a_2 \otimes \ldots \otimes a_n) \bigr),
\end{eqnarray*}}
showing that $\tilde{\phi}P_A=P\tilde{\phi}$. This equality will
be used in the sequel of the proof. We now show that
$\tilde{\phi}$ is an algebra morphism. This fact is verified for
words of length less or equal to 2. For words of higher lengths,
we set $X = x \otimes X'$ and $Y =y \otimes Y'$ and proceed by
induction on the total length of the considered words.
 \allowdisplaybreaks{
\begin{eqnarray*}
    \tilde{\phi}(a \otimes X \:\bar{\bullet}^{\ell}\: b \otimes Y)
    &=& \phi([a; b]) P\bigl( \tilde{\phi}( X \:\bullet^{\ell}\: Y) \bigr)\\
    &=& \phi([a; b]) P\bigl( \tilde{\phi}\bigl( x(X' \:\bullet^{\ell}\: Y) + y(X \:\bullet^{\ell}\: Y')
                                                    - [x; y] \otimes 1_A \otimes(X' \:\bullet^{\ell}\: Y')\bigr) \bigr)\\
    &=& \phi([a; b]) P\bigl( \tilde{\phi}(X \:\bar{\bullet}^{\ell}\: P_A(Y) + P_A(X) \:\bar{\bullet}^{\ell}\: Y
                                                    - X \:\bar{\bullet}^{\ell}\: (y \otimes 1_A \otimes Y')) \bigr)\\
    &=& \phi([a; b]) P\bigl( \tilde{\phi}(X)  \tilde{\phi}(P_A(Y))
                           + \tilde{\phi}(P_A(X)) \tilde{\phi}(Y)
                           - \tilde{\phi}(X) \tilde{\phi}(P_A(1_A))\tilde{\phi}(Y) \bigr)\\
    &=& \phi([a; b]) P\bigl( \tilde{\phi}(X)  P(\tilde{\phi}(Y))
                           + P(\tilde{\phi}(X)) \tilde{\phi}(Y)
                           - \tilde{\phi}(X) P(\tilde{\phi}(1_A))\tilde{\phi}(Y) \bigr)\\
    &=& \phi(a)\phi(b) P( \tilde{\phi}(X))  P( \tilde{\phi}(Y))
     = \tilde{\phi}(a \otimes X) \tilde{\phi}(b \otimes Y).
\end{eqnarray*}}

\noindent The $TD$-algebra morphism $\tilde{\phi}$ is unique since
it coincides with $\phi$ on $A$.
\end{proof}

\noindent There is a canonical functor \textbf{TDA}
$\longrightarrow$ \textbf{TriDend}, see \cite{Leroux04b}.

\begin{cor}  {\rm{\cite{Leroux04b}}} \label{cor:cdtFreeTD}
Let $A$ be a commutative algebra with unit $1_A$. The free unital
commutative $TD$-algebra over $A$, $\overline{\sha}^\ell(A)$, has
a canonical commutative tridendriform algebra structure given by,
$$
      X \prec Y := X \:\bar{\bullet}^{\ell}\: P_A(Y), \qquad
    X \bullet Y := - X \:\bar{\bullet}^{\ell}\: P_A(1_A) \:\bar{\bullet}^{\ell}\: Y.
$$
\end{cor}

\begin{rmk}{\rm{
If we drop the commutativity condition on $A$, then one still
finds a tridendriform algebra structure over the $TD$-algebra
$(T(A)^+, \bar{\bullet}^{\ell}, P_A)$. The third law is by
definition $X \succ Y:= P_A(X) \bar{\bullet}^{\ell} Y$.}}
\end{rmk}


\subsection{Commutative involutive \textbf{R}-algebras}
\label{ssect:invTDNA}

The aim of this part is to show that our results hold if one
considers the category of unital commutative algebras with
involution, denoted by \textbf{InvCom}. Recall that
\textbf{R}-algebra stands for either (scalar weight) Rota--Baxter,
Nijenhuis or $TD$-algebra.

\begin{thm} \label{thminv}
Let $A$ be a commutative algebra with unit $1_A$, equipped with an
involution $\dagger: A \rightarrow A$. Extending this involution
to $(T^+(A), \bar{\bullet}^{q}, P_A)$, for $q\in \{\theta,
\mathrm{r}, \ell \}$ by,
$$
    (a_1 \otimes \cdots \otimes a_n)^\dagger := a_1^\dagger \otimes \cdots \otimes a_n^\dagger,
$$
turns $(T^+(A), \bar{\bullet}^{q}, P_A, \dagger)$  into an
involutive commutative \textbf{R}-algebra which is also the free
involutive commutative \textbf{R}-algebra  over $A$.
\end{thm}

\begin{proof}
In fact, we only have to check the following. First, one observes
that $i_A([a; b]^\dagger) = i_A([b^\dagger; a^\dagger] =
i_A(b)^\dagger \:\bar{\bullet}^{\mathrm{q}}\: i_A(a)^\dagger =
i_A(a)^\dagger \:\bar{\bullet}^{\mathrm{q}}\: i_A(b)^\dagger$ as
expected. The \textbf{R}-map $P_A$ also behaves as expected, since
$P_A(U^\dagger)=1_A \otimes U^\dagger = ( 1_A \otimes U)^\dagger =
(P_A(U))^\dagger$. Let $(X,P,\dagger')$ be a unital commutative
\textbf{R}-algebra with involution. The \textbf{R}-algebra
morphism $\tilde{\phi}$ which is uniquely associated with the
algebra morphism $\phi: A \to X$ we found earlier, well-behaves
under the involution since $\tilde{\phi}(a^\dagger) = \phi
(a^\dagger) = \phi (a)^{\dagger'} = \tilde{\phi}(a)^{\dagger'},$
and by induction on the length of words, $\tilde{\phi}((a \otimes
U)^\dagger) = \tilde{\phi}(a^\dagger \otimes U^\dagger ) =
\tilde{\phi}(a^\dagger \:\bar{\bullet}^{\mathrm{q}}\:
P_A(U^\dagger)) =\tilde{\phi}(a)^{\dagger'} \:
\tilde{\phi}(P_A(U^\dagger)) = \tilde{\phi}(a)^{\dagger'}
P(\tilde{\phi}(U^\dagger)) = \tilde{\phi}(a)^{\dagger'}
P(\tilde{\phi}(U))^{\dagger'} = \tilde{\phi}(a \otimes
U)^{\dagger'},$ by using the commutativity of the product in $X$.
\end{proof}

\begin{cor} \label{cor:invCDT}
Let $(A, \dagger)$ be a commutative involutive algebra with unit
$1_A$. Extend $\dagger$ as explained in Theorem~\ref{thminv}.
Then, $(T^+(A), \bar{\bullet}^{\ell}, P_A, \dagger)$ has an
involutive commutative tridendriform algebra structure, i.e., the
two binary operations given by,
\begin{equation}
    \label{eqs:invCDT}
    X \prec Y := X \:\bar{\bullet}^{\ell}\: P_A(Y),
    \qquad
    X \bullet Y := - X \: \bar{\bullet}^{\ell}\: P_A(1_A) \: \bar{\bullet}^{\ell}\: Y,
\end{equation}
verify $ (X \prec Y)^\dagger =  X^\dagger \prec Y^\dagger$ and $(X
\bullet Y)^\dagger=  X^\dagger \bullet Y^\dagger$.
\end{cor}


\section{Quasi-shuffle product from left-shift shuffle product}
\label{sect:QshufLshuf}

We present here a relation between the quasi-shuffle and the
left-shift shuffle. Let $V$ be a $\mathbb{K}$-vector space. Recall
that $T(V)$ is the free unital associative $\mathbb{K}$-algebra
over $V$ with unit $1_\mathbb{K}$, and
$$
    S(V):= \mathbb{K} \cdot 1_\mathbb{K} \oplus V \oplus \bigoplus_{n>1} S^n(V),
$$
is the free unital commutative algebra over $V$, where $S^n(V)$ is
the quotient of $V^{\otimes n}$ by the action of the symmetric
group $S_n$. Set
$$
    \bar{S}(V) := S(V)/\mathbb{K} \cdot 1_\mathbb{K}= V \oplus \bigoplus_{n>1} S^n(V)
     \; \makebox{ and }\;
    \bar{T}(V):=T(V)/\mathbb{K} \cdot 1_\mathbb{K}=V \oplus \bigoplus_{n>1} V^{\otimes n}.
$$

\noindent By $F(V)$ ($\bar{F}(V)$), we mean either $S(V)$ or
$T(V)$ (resp. $\bar{S}(V)$ or $\bar{T}(V)$) and by $F^n(V)$
($\bar{F}^n(V)$) either $S^n(V)$ or $V^{\otimes n}$, $n\geq 0$
(resp. $n>0$).

In this section we fix the quasi-shuffle weight to be $\theta=1$,
see Eq.~(\ref{eq:geneShuffle}). By $\bar{T}_1(\bar{S}(V))$, (resp.
$\bar{T}_1(\bar{T}(V))$), we mean the vector-space
$\bar{T}(\bar{S}(V))$, (resp. $\bar{T}(\bar{T}(V))$) equipped with
binary operations $\prec_1$, $\bullet_1$, (resp. with $\succ_1$
being added) defined in terms of the quasi-shuffle product of
weight one,
$$
           a\otimes X \prec_1 b\otimes Y  := a \otimes (X \bullet^1 b\otimes Y),
    \quad a\otimes X \bullet_1 b\otimes Y := [a;b] \otimes (X \bullet^1 Y),
$$
$$
   (\textrm{and} \ a\otimes X \succ_1 b\otimes Y := b \otimes (a\otimes X \bullet^1 Y) \ \textrm{to\ be\ added}),
$$

\noindent recall that $\bullet^1$ stands for the quasi-shuffle
product of weight one, see relation~(\ref{eq:mixShuf}), and that
$T^+(F(V)):=F(V) \otimes T(F(V))$. Here, $[-;-]$ is the product in
$F(V)$ which we will denote by concatenation, $[a;b]=ab$, $a,b \in
F(V)$.\smallskip

Define the linear map $\Omega: \bar{F}^n(V) \rightarrow{}
\bar{F}^n(V) \otimes \mathbb{K} \cdot 1_\mathbb{K}^{\otimes
(n-1)}$, by $v_1 \cdots v_n \mapsto (-1)^{n+1} v_1 \cdots v_n
\otimes 1_\mathbb{K}^{\otimes (n-1)}$. Extend this map as follows,
$$
    \Omega: \bar{T}_1(\bar{F}(V)) \rightarrow (T^+(F(V)), \prec,\bullet,\succ)
$$
$$
    a_{1} \otimes \ldots \otimes a_{k} \mapsto \Omega(a_{1}) \otimes \ldots \otimes \Omega(a_{k}).
$$
with $\prec,\bullet$ and $\succ$ defined in terms of the
left-shuffle product $\bar{\bullet}^\ell$ and the right-shift map
$P_{F(V)}$
\begin{equation*}
    X \prec Y := X \:\bar{\bullet}^{\ell}\: P_{F(V)}(Y),
    \qquad
    X \succ Y := P_{F(V)}(X) \:\bar{\bullet}^{\ell}\: Y,
    \qquad
    X \bullet Y := - X \: \bar{\bullet}^{\ell}\: P_{F(V)}(1_{\mathbb{K}}) \: \bar{\bullet}^{\ell}\: Y.
\end{equation*}

\noindent Here, the letters $a_{i}$ are words in $\bar{F}(V)$,
i.e., $a_{i}=v_{i,1} \cdots v_{i,m_i} \in V^{\otimes m_i}$,
$m_i>0$. We use the word notation for elements in $\bar{F}(V)$,
that is, we omit the tensor product sign, which we reserve for
either $\bar{T}_1(\bar{F}(V))$ or $T^+(F(V))$.

\begin{prop} \label{Loday1}
The linear map $\Omega$, so defined, is an injective morphism for
the two (resp. three) operations, $\prec_1, \bullet_1$ (resp.
$\succ_1$) just defined.
\end{prop}

\begin{proof} We restrict ourselves to the operations $\prec_1,
\bullet_1$. Injectivity is readily verified. We proceed by
induction.
 \allowdisplaybreaks{
\begin{eqnarray*}
    \Omega(v_1 \prec_1 v_2) = \Omega(v_1 \otimes v_2) &=& \Omega(v_1) \otimes \Omega(v_2)\\
                                            &=& v_1 \otimes v_2= v_1 \:\bar{\bullet}^{\ell}\: 1_{\mathbb{K}} \otimes v_2\\
                                            &=& v_1 \:\bar{\bullet}^{\ell}\: P_{F(V)}(v_2) \\
                                            &=& \Omega(v_1) \prec \Omega(v_2),
\end{eqnarray*}}
since $\Omega$ is the identity on $V$. Similarly,
 \allowdisplaybreaks{
\begin{eqnarray*}
    \Omega(v_1 \bullet_1 v_2) &=& \Omega([v_1;v_2]) =-[v_1;v_2] \otimes 1_\mathbb{K}\\
                         &=& -v_1v_2 \otimes 1_\mathbb{K}\\
                         &=& -v_1 \:\bar{\bullet}^{\ell}\: (1_{\mathbb{K}} \otimes 1_{\mathbb{K}})
                                  \:\bar{\bullet}^{\ell}\: v_2
                          =  v_1 \bullet v_2 \\
                         &=& \Omega(v_1) \bullet \Omega(v_2).
\end{eqnarray*}}

\noindent Recall that the product $[-;-]$ here is the one in
$T(V)$ ($S(V)$). Now,
 \allowdisplaybreaks{
\begin{eqnarray*}
    \Omega(a\otimes X \prec_1 b\otimes Y) &=& \Omega \bigl( a \otimes (X \:\bullet^1\: (b \otimes Y) ) \bigr)\\
                                     &=& \Omega(a) \otimes \Omega(X \:\bullet^1\: (b \otimes Y))
                                      = \Big((-1)^{l(a)+1} a \otimes 1_\mathbb{K}^{\otimes (l(a)-1)}\Big)
                                                                       \otimes \Omega( X \:\bullet^1\: (b\otimes Y)),\\
        &=& (-1)^{l(a)+1} a \otimes 1_\mathbb{K}^{\otimes (l(a)-1)} \otimes \Omega(X \:\bullet^1\: b \otimes Y)\\
        &=& (-1)^{l(a)+1} a \otimes \bigl(1_\mathbb{K}^{\otimes (l(a)-1)} \:\bullet^\ell\: \Omega(X \:\bullet^1\: b \otimes Y)\bigr)\\
        &=& (-1)^{l(a)+1} a \otimes \bigl(1_\mathbb{K}^{\otimes (l(a)-1)} \:\bullet^\ell\: \Omega(X) \:\bullet^\ell\: \Omega(b \otimes Y)\bigr)\\
        &=& (-1)^{l(a)+1} \bigl(a \otimes \bigl(1_\mathbb{K}^{\otimes (l(a)-1)} \:\bullet^\ell\: \Omega(X)\bigr)
                                                              \:\bar{\bullet}^\ell\: P_{F(V)}\bigl(\Omega(b \otimes Y)\bigr),\\
        &=& (-1)^{l(a)+1} \Big( a \otimes \bigl(1_\mathbb{K}^{\otimes (l(a)-1)} \:\bullet^\ell\: \Omega(X) \bigr)\Big) \prec \Omega(b \otimes Y),\\
        &=& (-1)^{l(a)+1} \Big(\bigl( a \otimes 1_\mathbb{K}^{\otimes (l(a)-1)} \bigr)  \:\bar{\bullet}^\ell\: P_{F(V)}\bigl(\Omega(X)\bigr)\Big)
                                                                     \prec \Omega(b \otimes Y),\\
        &=& (-1)^{l(a)+1} \Big(\bigl( a \otimes 1_\mathbb{K}^{\otimes (l(a)-1)} \bigr) \prec \Omega(X) \Big) \prec \Omega(b \otimes Y),\\
        &=& (\Omega(a) \prec \Omega(X)) \prec \Omega(b \otimes Y) = \Omega(a \otimes X) \prec \Omega(b \otimes Y).
\end{eqnarray*}}

\noindent We used Eqs.~(\ref{eq:lID}), and the induction
hypothesis in the part $\Omega\bigl(X \:\bullet^1\: (b \otimes
Y)\bigr)= \Omega(X) \:\bullet^\ell\: \Omega(b \otimes Y)$, since
$\bullet^1$ and $\bullet^\ell$ are the sum of the dendriform
operations $\prec_1, \bullet_1$ and $\prec, \bullet$,
respectively. For the $\bullet$ operation, we find on the one
hand,
 \allowdisplaybreaks{
\begin{eqnarray*}
    \Omega(a\otimes X \bullet_1 b\otimes Y) &=& \Omega\bigl([a;b] \otimes (X \bullet^1 Y)\bigr)
                                        =  \Omega([a;b]) \otimes \Omega(X \bullet^1 Y)\\
                                       &=& \Omega([a;b]) \otimes \bigl(\Omega(X) \:\bullet^\ell\: \Omega(Y)\bigr),\\
                                       &=& (-1)^{(l(a)+l(b)+1)} [a;b] \otimes 1_{\mathbb{K}}^{\otimes (l(a)+l(b)-1)}
                                                                      \otimes \bigl(\Omega(X)\bullet^\ell \Omega(Y)\bigr).
\end{eqnarray*}}

\noindent On the other hand,
 \allowdisplaybreaks{
\begin{eqnarray*}
 \Omega(a\otimes X) \bullet \Omega(b\otimes Y) &=& \bigl(\Omega(a)\otimes \Omega(X)\bigr) \bullet \bigl(\Omega(b)\otimes \Omega(Y)\bigr)
                                      =  -\bigl(\Omega(a)\otimes \Omega(X)\bigr)\:\bar{\bullet}^\ell\: \bigl(1_{\mathbb{K}}^{\otimes 2}\bigr)
                                                              \:\bar{\bullet}^\ell\: \bigl(\Omega(b)\otimes \Omega(Y)\bigr),\\
                                     &=& -(-1)^{ (l(a)+l(b)+1+1)} \bigl(a \otimes 1_{\mathbb{K}}^{\otimes (l(a)-1)} \otimes \Omega(X)\bigr)
                                               \:\bar{\bullet}^\ell\: \bigl(1_{\mathbb{K}}^{\otimes 2}\bigr)
                                               \:\bar{\bullet}^\ell\:
                                                  \bigl(b \otimes 1_{\mathbb{K}}^{\otimes (l(b)-1)} \otimes \Omega(Y)\bigr),\\
                                     &=& (-1)^{ (l(a)+l(b)+1)} [a;b]\otimes 1_{\mathbb{K}}^{\otimes (l(a)-1)+(l(b)-1)+1}
                                                              \otimes (\Omega(X) \:\bullet^\ell\: \Omega(Y)),
\end{eqnarray*}}
hence the equality. Similarly, we prove $\Omega(a\otimes X \succ_1
b\otimes Y) = \Omega(a\otimes X)\succ  \Omega(b\otimes Y)$.
\end{proof}

\begin{cor} \label{Loday2}
Let $V$ be a $\mathbb{K}$-vector space. Then, any monomial of
$\bar{T}_1(\bar{S}(V))$
          (resp. $\bar{T}_1(\bar{T}(V))$), $a_1 \otimes \ldots \otimes a_n$, can be written as,
$$
    a_1 \prec_1  (a_2 \prec_1 ( \ldots (a_{n-1} \prec_1 a_n)) \cdots),
$$
where $a_i:=v_{i,1} \otimes \ldots \otimes v_{i,k_i}$ can be
written as $a_i =v_{i,1} \bullet_1 \ldots \bullet_1 v_{i,k_i}$.
\end{cor}

\begin{proof}
Let $a_1 \otimes \ldots \otimes a_n$ be such a monomial. Observe
that $\Omega(a_1 \otimes \ldots \otimes a_n)$ can always be
written as
$$
    \Omega(a_1) \prec ( \Omega(a_2) \prec  \ldots (\Omega(a_{n-1}) \prec  \Omega(a_n)) \cdots),
$$
where $\Omega(a_i):=(-1)^{(k_i + 1)}(v_{i,1} \otimes \ldots
\otimes v_{i,k_i})\otimes 1_{\mathbb{K}}^{\otimes (k_i - 1)}$, can
be rewritten as $v_{i,1} \bullet \ldots \bullet v_{i,k_i}$. Apply
now the inverse of $\Omega$ to conclude.
\end{proof}

\noindent The fact that $\bar{T}_q(\bar{F}(V))$ has a
(commutative) tridendriform algebra structure was found by Loday
and Ronco in~\cite{LodayRonco04}. The relation in item two of
Corollary~\ref{Loday2} was proven inductively by Loday
in~(\cite{Loday05}). Here, we related Loday's result to the
left-shift shuffle product.\\

We finish this section by presenting a different proof of a
theorem of Loday \cite{Loday05}. The free C$TD$-algebra (resp. the
free involutive tridendriform algebra) over a $\mathbb{K}$-vector
space $V$ is a C$TD$-algebra $CTD(V)$, (resp. an involutive
tridendriform algebra $ITD(V)$) equipped with a linear map $i_V:V
\rightarrow CTD(V)$ (resp. $i_V:V \rightarrow ITD(V)$) such that
for any linear map $\phi: V \rightarrow A$, where $A$ is a
C$TD$-algebra, (resp. an involutive tridendriform algebra), there
exists a unique C$TD$-algebra morphism, (resp. an involutive
tridendriform algebra morphism) $\tilde{\phi}: CTD(V) \rightarrow
A$ (resp. $\tilde{\phi}: ITD(V) \rightarrow A$) such that
$\tilde{\phi} \circ i_V= \phi$. In the involutive case, we remind
the reader that the involution on $\bar{T}(V)$ is defined to be
$(v_1 \dots v_n)^\dagger:= v_n \dots v_1$, we omitted tensor
signs.

\begin{thm}\cite{Loday05}
Let $V$ be a $\mathbb{K}$-vector space, then
$\bar{T}_1(\bar{S}(V))$ (resp. $(\bar{T}_1(\bar{T}(V)), \dagger)$)
is the free C$TD$-algebra (resp. involutive tridendriform algebra)
over $V$.
\end{thm}

\begin{proof}
Let $A$ be a C$TD$-algebra with operations ($\prec$, $\bullet$)
(resp. involutive tridendriform algebra with operations ($\prec$,
$\succ$, $\bullet$) and involution $\dagger$) and $\phi: V
\rightarrow A$ be a linear map. We denote the product
(\ref{def:double}) in $A$ by $a*b:=a\prec b + a \succ b + a
\bullet b$.

Recall that in $\bar{T}_1(\bar{F}(V))$ we may write $a_1 \otimes
\ldots \otimes a_{n-1} \otimes a_n = a_1 \prec_1 \bigl( \ldots
\prec_1 (a_{n-1} \prec_1 a_n) \cdots \bigr)$, where $a_i =v_{i,1}
\bullet_1 \ldots \bullet_1 v_{i,k_i}$ for $a_i=v_{i,1} \cdots
v_{i,k_i} \in \bar{F}(V)$. Set
$$
    \tilde{\phi}(a_1 \otimes \ldots \otimes a_{n-1} \otimes a_n) :=
    \tilde{\phi}(a_1) \prec \bigl( \tilde{\phi}(a_2) \prec\ldots (\tilde{\phi}(a_{n-1}) \prec \tilde{\phi}(a_n) ) \cdots\bigr),
$$
where $\tilde{\phi}(a_i):=\phi(v_{i,1}) \bullet \ldots \bullet
\phi(v_{i,k_i})$. So defined, $\tilde{\phi}$ is a morphism of
C$TD$-algebras (resp. involutive tridendriform algebras). We
proceed by induction on the total length of words. The result is
trivial for a total length of two. Therefore,
 \allowdisplaybreaks{
\begin{eqnarray*}
    \tilde{\phi}\bigl(a \otimes X \prec_1 (b\otimes Y)\bigr) &=& \tilde{\phi}\bigl(a\otimes (X \bullet^1 (b\otimes Y))\bigr)\\
                                                      &=&  \tilde{\phi}(a) \prec \tilde{\phi}\bigl(X \bullet^1 (b \otimes Y)\bigr)\\
                                                      &=& \tilde{\phi}(a) \prec \bigl(\tilde{\phi}(X) * \tilde{\phi}(b \otimes Y)\bigr)\\
                                                      &=& \bigl(\tilde{\phi}(a) \prec \tilde{\phi}(X)\bigr) \prec \tilde{\phi}(b\otimes Y)\\
                                                      &=& \tilde{\phi}(a\otimes X) \prec \tilde{\phi}(b\otimes Y).\\
\end{eqnarray*}}
In the case of involutive tridendriform algebras, we obtain in
addition the equality $\tilde{\phi}(a\otimes X \succ_1 b\otimes
Y)= \tilde{\phi}(a\otimes X) \succ_1 \tilde{\phi}(b\otimes Y)$. As
for the operation $\prec_1$ we find for $\bullet_1$,
 \allowdisplaybreaks{
\begin{eqnarray*}
    \tilde{\phi}(a\otimes X \bullet_1 b\otimes Y) &=& \tilde{\phi}\bigl([a;b]\otimes (X \bullet^1 Y)\bigr)
                                                   = \tilde{\phi}([a;b]) \prec \tilde{\phi}(X \bullet^1 Y)
                                                   = \tilde{\phi}([a;b]) \prec \bigl(\tilde{\phi}(X) * \tilde{\phi}(Y)\bigr) \\
    &=& \bigl(\tilde{\phi}(a) \bullet \tilde{\phi}(b)\bigr) \prec \bigl(\tilde{\phi}(X) * \tilde{\phi}(Y)\bigr)
     = \Big( \bigl(\tilde{\phi}(a) \bullet \tilde{\phi}(b)\bigr) \prec \tilde{\phi}(X)\Big) \prec \tilde{\phi}(Y).
\end{eqnarray*}}
We used the induction hypothesis in the part $\tilde{\phi}(X
\bullet^1 Y) = \tilde{\phi}(X) * \tilde{\phi}(Y)$. Now, using
commutativity, we find
 \allowdisplaybreaks{
\begin{eqnarray*}
    &=& \bigl((\tilde{\phi}(b) \bullet \tilde{\phi}(a)) \prec \tilde{\phi}(X)\bigr) \prec \tilde{\phi}(Y))
     =  \bigl( \tilde{\phi}(b) \bullet ( \tilde{\phi}(a) \prec \tilde{\phi}(X) ) \bigr) \prec \tilde{\phi}(Y)\\
    &=& \bigl(\tilde{\phi}(b) \bullet \tilde{\phi}(a \otimes X)\bigr) \prec \tilde{\phi}(Y)
     =  \bigl(\tilde{\phi}(a\otimes X) \bullet \tilde{\phi}(b)\bigr) \prec \tilde{\phi}(Y)
     =  \tilde{\phi}(a\otimes X) \bullet (\tilde{\phi}(b) \prec \tilde{\phi}(Y))\\
    &=& \tilde{\phi}(a\otimes X) \bullet \tilde{\phi}(b\otimes Y).
\end{eqnarray*}}
Now, in the involutive case, we compute
 \allowdisplaybreaks{
\begin{eqnarray*}
    \tilde{\phi}(a\otimes X) \bullet \tilde{\phi}(b\otimes Y) &=&
                        \bigl(\tilde{\phi}(a\otimes X) \bullet \tilde{\phi}(b\otimes Y)\bigr)^{\dagger\dagger}
                         = \Big( \bigl(\tilde{\phi}(a) \prec \tilde{\phi}(X)\bigr) \bullet \bigl(\tilde{\phi}(b\otimes Y)\bigr)\Big)^{\dagger\dagger}\\
    &=& \big[ \tilde{\phi}(b\otimes Y)^{\dagger} \bullet \bigl(\tilde{\phi}(X)^{\dagger} \succ \tilde{\phi}(a) ^{\dagger}\bigr)\big]^\dagger
     =  \big[\bigl(\tilde{\phi}(b\otimes Y)^{\dagger} \prec \tilde{\phi}(X)^{\dagger}\bigr) \bullet \tilde{\phi}(a)^{\dagger}\big]^\dagger \\
    &=& \big[\bigl((\tilde{\phi}(b) \prec \tilde{\phi}(Y))^{\dagger} \prec \tilde{\phi}(X)^{\dagger}\bigr)
                                        \bullet \tilde{\phi}(a)^{\dagger}\big]^\dagger \\
    &=& \big[\bigl((\tilde{\phi}(Y)^{\dagger} \succ \tilde{\phi}(b)^{\dagger}) \prec \tilde{\phi}(X)^{\dagger}\bigr)
                                        \bullet \tilde{\phi}(a)^{\dagger}\big]^\dagger \\
    &=& \big[\bigl(\tilde{\phi}(Y)^{\dagger} \succ ( \tilde{\phi}(b)^{\dagger} \prec \tilde{\phi}(X)^{\dagger})\bigr)
                                        \bullet \tilde{\phi}(a)^{\dagger}\big]^\dagger \\
    &=& \big[\bigl(\tilde{\phi}(Y)^{\dagger} \succ \bigl( (\tilde{\phi}(X) \succ \tilde{\phi}(b))^{\dagger}\bigr)\bigr)
                                        \bullet \tilde{\phi}(a)^{\dagger}\big]^\dagger\\
    &=& \tilde{\phi}(a) \bullet \bigl(( \tilde{\phi}(b) \prec \tilde{\phi}(X)) \prec \tilde{\phi}(Y)\bigr)
     = \bigl(\tilde{\phi}(a) \bullet \bigl(\tilde{\phi}(b) \prec \tilde{\phi}(X)\bigr)\bigr) \prec \tilde{\phi}(Y)\\
    &=& \bigl((\tilde{\phi}(a) \bullet  \tilde{\phi}(b)) \prec \tilde{\phi}(X)\bigr) \prec \tilde{\phi}(Y).
\end{eqnarray*}}
Hence the expected equality. The uniqueness of $\tilde{\phi}$ is
readily verified, which finishes the proof.
\end{proof}


\section{$\textrm{{\bf{R}}}$-bialgebras}
\label{sect:Rbialg}

In Subsection \ref{ssect:RBcategories}, and Paragraphs
\ref{sssect:freeRBandNA} and \ref{sssect:freeTD}, the categories
of $\textrm{{\bf{R}}}$-algebras were defined, and their free
objects constructed, respectively.

The aim of this part is to address the following question. If we
start with a commutative bialgebra $A$ instead of an algebra, how
much of its bialgebra structure can be lifted to $T^+(A)$ equipped
with one of the generalized shuffle products? In fact we will see
that this motivates the introduction of the notion of right (left)
$\textrm{{\bf{R}}}$-bialgebras.

Recall that $(\mathbb{K}, -\theta \id_{\mathbb{K}})$ is in
\textbf{ComRBA} and $(\mathbb{K}, \id_{\mathbb{K}})$ is naturally
in \textbf{ComNA} and \textbf{ComTDA}. Next, observe that if $(B,
P)$ is an $\textrm{{\bf{R}}}$-algebra then for each $n>0$,
$B^{\otimes n}$, equipped with the linear map $P^{(n)} := P
\otimes \id_B^{\otimes (n-1)}: B^{\otimes n} \rightarrow
B^{\otimes n}$ turns into an $\textrm{{\bf{R}}}$-algebra. One can
also choose to define $\hat{P}^{(n)}:= \id_B^{\otimes (n-1)}
\otimes P : B^{\otimes n} \rightarrow B^{\otimes n}$.

\begin{defn} \label{def:Rbialgebra}
The $\textrm{{\bf{R}}}$-algebra $(B, P)$ is said to be a right
$\textrm{{\bf{R}}}$-bialgebra if it is equipped with a
coassociative coproduct $\Delta: B \rightarrow B^{\otimes 2}$ and
a right counit $\epsilon : B \rightarrow \mathbb{K}$, i.e.,
verifying $(\id_B \otimes \epsilon) \Delta = \id_B$, which are
both $\textrm{{\bf{R}}}$-algebra morphisms.
\end{defn}

\noindent  In case we defined the map $\hat{P}^{(n)}:=
\id_B^{\otimes (n-1)} \otimes P$, we obtain a left
$\textrm{{\bf{R}}}$-bialgebra, i.e., with a left counit instead of
a right one.


\subsection{Case 1: bialgebras}
\label{ssect:case1}

Recall that
$$
    T^+(H) := H \otimes T(H) = H \otimes \mathbb{K} 1_\mathbb{K} \oplus \bigoplus_{n>1} H^{\otimes n},
$$
and that we denote by $(T^+(H),\bar{\bullet}^q,P_H)$, with $q \in
\{\theta, {\mathrm{r}},\ell \}$, the free commutative
$\textrm{{\bf{R}}}$-algebra over a unital commutative algebra $H$,
and $P_H$ stands for the right-shift operator.

\begin{thm} \label{thmHopfnij}
Let $(H,m_H,\delta,\epsilon)$ be a commutative bialgebra with unit
$1_H$. Then, the free commutative $\textrm{{\bf{R}}}$-algebra
$(T^+(H),\bar{\bullet}^q,P_H)$ has the structure of a right
$\textrm{{\bf{R}}}$-bialgebra.
\end{thm}

\begin{proof}
Recall that $(T^+(H)^{\otimes 2},\bar{\bullet}^q,P^{(2)}_H)$ is
also an $\textrm{{\bf{R}}}$-algebra. Here we keep
$\bar{\bullet}^q$ to denote the product in $T^+(H)^{\otimes 2}$.
Let $i_H: H \hookrightarrow (T^+(H),\bar{\bullet}^q,P_H)$ be the
inclusion map, defined naturally as
\begin{equation}
\label{inclusionmap}
    a \mapsto i_H(a) := a \otimes 1_{\mathbb{K}} \in T^+(H),
\end{equation}
which is an algebra morphism, see~(\ref{inclusion}). Consider the
following commutative diagram
\begin{diagram}
    H          & \rInto^{\phantom{mim} i_H \phantom{mim}}  & (T^+(H),\bar{\bullet}^q,P_H)\\
               & \rdTo_{\delta\ }                        & \dTo_{ \Delta\ }\\
               &                                         & (T^+(H)^{\otimes 2},\bar{\bullet}^q,P^{(2)}_H)
\end{diagram}
Here we understand the map $\delta$, keeping the same notation, as
derived from the coproduct on $H$ as follows
$$
    \delta : H \to T^+(H)^{\otimes 2},\quad a \mapsto \delta(a):=
    \bigl(a_{(1)} \otimes 1_{\mathbb{K}}\bigr) \otimes \bigl(a_{(2)} \otimes 1_{\mathbb{K}}\bigr),
$$
using Sweedler's notation for $\delta(a)= a_{(1)} \otimes a_{(2)}$
on $H$. So defined, $\delta$ is a unital algebra morphism. Indeed,
we find immediately that
 \allowdisplaybreaks{
\begin{eqnarray*}
    \delta\bigl(m_H(a\otimes b)\bigr) &=& \bigl( i_H(m_H(a_{(1)}\otimes b_{(1)}))\bigr)
                                        \otimes
                                        \bigl(i_H(m_H(a_{(2)}\otimes b_{(2)}))\bigr)\\
                                    &=& \delta(a) \:\bar{\bullet}^q\: \delta(b)
\end{eqnarray*}}
and hence there exists a unique $\textrm{{\bf{R}}}$-algebra map
$\Delta$ which extends $\delta$, such that $\Delta \circ i_H
=\delta$. Furthermore, observe that $(T^+(H)^{\otimes
3},\bar{\bullet}^q,P_H^{(3)})$ also is a commutative
$\textrm{{\bf{R}}}$-algebra. Now, let us define the following map
$$
 \delta^{(2)}:= i_H^{\otimes 3} \circ (\id_H \otimes \delta)\delta: H \to T^+(H)^{\otimes 3},
$$
which is a unital algebra morphism. Due to the coassociativity of
$\delta$ on $H$ we have that $\delta^{(2)}= i_H^{\otimes 3} \circ
(\delta \otimes \id_H)\delta=:\hat{\delta}^{(2)}$. Define
$\Delta^{(2)}:=(\id_{T^+(H)} \otimes \Delta)\Delta$ and
$\hat{\Delta}^{(2)}:=(\Delta \otimes \id_{T^+(H)})\Delta$ and
consider the following commutative diagram
\begin{diagram}
    H          & \rInto^{\phantom{mim} i_H \phantom{mim}}  & (T^+(H),\bar{\bullet}^q,P_H)\\
               & \rdTo_{\delta^{(2)}\ }                    & \dTo_{\Delta^{(2)}}\\
               &                                           & (T^+(H)^{\otimes 3},\bar{\bullet}^q,P_H^{(3)})
\end{diagram}
such that both $\delta^{(2)}$ and $\hat{\delta}^{(2)}$ extend to
the $\textrm{{\bf{R}}}$-algebra maps $(\id_{T^+(H)} \otimes
\Delta) \Delta$ and $(\Delta \otimes \id_{T^+(H)}) \Delta$,
respectively. Hence, coassociativity of $\Delta$ follows by
unicity of the extension. Consequently, the map $\Delta$ is both a
coassociative coproduct and an $\textrm{{\bf{R}}}$-algebra
morphism in the sense of Definition~\ref{def:Rbialgebra}.

As $(\mathbb{K}, -\theta \id_\mathbb{K})$, respectively
$(\mathbb{K}, \id_\mathbb{K})$ are commutative
$\textrm{{\bf{R}}}$-algebras, the unital morphism $\epsilon$ has a
unique $\textrm{{\bf{R}}}$-algebra map extension still denoted by
$\epsilon$ such that $\epsilon \circ i_H=\epsilon$.

As $i_H: H \rightarrow ({T^+(H)},\bar{\bullet}^q,P_H)$ is a unital
algebra morphism, it has a unique $\textrm{{\bf{R}}}$-algebra map
extension which is, of course, $\id_{T^+(H)} :
(T^+(H),\bar{\bullet}^q,P_H) \rightarrow
(T^+(H),\bar{\bullet}^q,P_H)$. Now observe that $(\id_{T^+(H)}
\otimes \epsilon) \Delta: (T^+(H),\bar{\bullet}^q,P_H) \rightarrow
(T^+(H),\bar{\bullet}^q,P_H)$ is an $\textrm{{\bf{R}}}$-algebra
map since
$$
    (\id_{T^+(H)}  \otimes \epsilon)\Delta \circ P_H = (\id_{T^+(H)}  \otimes \epsilon)(P_H \otimes \id_{T^+(H)} ) \Delta
                                                     = P_H \circ (\id_{T^+(H)}  \otimes \epsilon) \Delta.
$$
Moreover, observe that
$$
    (\id_{T^+(H)} \otimes \epsilon)\Delta \circ i_H = (\id_{T^+(H)} \otimes \epsilon)\delta
                                                    = i_H.
$$
Consequently, $(\id_{T^+(H)} \otimes \epsilon)\Delta
=\id_{T^+(H)}$ by unicity of the extension of the inclusion $i_H$.

However, observe that $(\epsilon \otimes \id_{T^+(H)})\Delta:
(T^+(H),\bar{\bullet}^q,P_H) \rightarrow
(T^+(H),\bar{\bullet}^q,P_H)$ is an algebra morphism but not an
$\textrm{{\bf{R}}}$-algebra map since
$$
    (\epsilon \otimes \id_{T^+(H)})\Delta \circ P_H = (\epsilon \otimes \id_{T^+(H)})(P_H \otimes \id_{T^+(H)}) \Delta
                                                    = (\id_{\mathbb{K}} \circ \epsilon \otimes \id_{T^+(H)}) \Delta
                                                    = (\epsilon \otimes \id_{T^+(H)}) \Delta.
$$
Therefore, we get a right $\textrm{{\bf{R}}}$-bialgebra structure
on $(T^+(H),\bar{\bullet}^q,P_H)$.
\end{proof}

\begin{rmk} {\rm{
\begin{itemize}

    \item{ Similarly, if we work with a commutative Hopf algebra,
the antipode $S: H \rightarrow H$, which is a unital algebra
morphism, can be uniquely lifted to an $\textrm{{\bf{R}}}$-algebra
map $\tilde{S}: (T^+(H),\bar{\bullet}^q,P_H) \rightarrow
(T^+(H),\bar{\bullet}^q,P_H)$ verifying $\tilde{S} \circ i_H = S$.
However the two defining antipode equations, i.e.,
$\bar{\bullet}^q(\id_{T^+(H)} \otimes \tilde{S})\Delta = \eta
\circ \epsilon = \bar{\bullet}^q(\tilde{S} \otimes
\id_{T^+(H)})\Delta$, where $\eta: \mathbb{K} \rightarrow
(T^+(H),\bar{\bullet}^q,P_H)$, $k \mapsto k \cdot 1_H\otimes
1_\mathbb{K}$ is the unit map, are lost.}

    \item{Observe that we considered here only the categories of
commutative Rota--Baxter, Nijenhuis and $TD$-algebras, but, the
proof does not use any of those properties. In fact, the result of
this theorem is still valid for any category of commutative
algebras with objects $(B,X)$, where $X: B \rightarrow B$ is a
linear map, and morphisms as in
Subsection~\ref{ssect:RBcategories}. In fact the only constraint
is that $\mathbb{K}$ is of $X$-type and that there exists a
universal object in this category.}

\end{itemize}}}
\end{rmk}

\noindent By an $\textrm{{\bf{R}}}$-bialgebra morphism $f : B_1
\rightarrow B_2$, we mean an $\textrm{{\bf{R}}}$-algebra morphism
which is also a coalgebra morphism. The category of commutative
right $\textrm{{\bf{R}}}$-bialgebras is well-defined. Theorem
\ref{thmHopfnij} can be rephrased by the following corollary.

\begin{cor}
Let $q \in \{\theta,{\mathrm{r}},\ell\}$. The covariant functor
$\overline{\sha}^\mathrm{q}$ embeds the category of commutative
bialgebras into the category of commutative right
$\textrm{{\bf{R}}}$-bialgebras.
\end{cor}

\begin{proof}
Let $f: H_1 \rightarrow H_2$ be a morphism between two commutative
bialgebras. We know that $\overline{\sha}^\mathrm{q}(H_1)$ and
$\overline{\sha}^\mathrm{q}(H_2)$ are commutative right
$\textrm{{\bf{R}}}$-bialgebras. With notation of Theorem
\ref{thmHopfnij}, the unital algebra morphism $i_2 \circ f: H_1
\rightarrow \overline{\sha}^\mathrm{q}(H_2)$ can be lifted into a
unique $\textrm{{\bf{R}}}$-algebra morphism $\zeta:
\overline{\sha}^\mathrm{q}(H_1) \rightarrow
\overline{\sha}^\mathrm{q}(H_2)$ verifying $\zeta \circ i_1 = i_2
\circ f$. We rename $\zeta$ into $\overline{\sha}^\mathrm{q}(f)$.
Now, $(f \otimes f) \delta_1:H_1 \rightarrow
\overline{\sha}^\mathrm{q}(H_2)^{\otimes 2}$ is also a unital
algebra morphism which can be lifted into a unique
$\textrm{{\bf{R}}}$-algebra morphism $\Gamma:
\overline{\sha}^\mathrm{q}(H_1) \rightarrow
\overline{\sha}^\mathrm{q}(H_2)^{\otimes 2}$ verifying $\Gamma
\circ i_1 = (f \otimes f) \delta_1$. Similarly for the unital
algebra map $\delta_2\circ f:H_1 \rightarrow
\overline{\sha}^\mathrm{q}(H_2)^{\otimes 2}$, we get a unique
$\textrm{{\bf{R}}}$-algebra morphism $\Gamma':
\overline{\sha}^\mathrm{q}(H_1) \rightarrow
\overline{\sha}^\mathrm{q}(H_2)^{\otimes 2}$ verifying $\Gamma'
\circ i_1 = \delta_2\circ f$. Thus by unicity $\Gamma=\Gamma'$
since $f$ is a coalgebra morphism. But the
$\textrm{{\bf{R}}}$-algebra morphism $\Delta_2
\overline{\sha}^\mathrm{q}(f)$ verifies $\Delta_2
\overline{\sha}^\mathrm{q}(f) \circ i_1 =\Delta_2 i_2 \circ f =
\delta_2 \circ f $ and the $\textrm{{\bf{R}}}$-algebra morphism
$(\overline{\sha}^\mathrm{q}(f) \otimes
\overline{\sha}^\mathrm{q}(f))\Delta_1$ verifies
$(\overline{\sha}^\mathrm{q}(f) \otimes
\overline{\sha}^\mathrm{q}(f))\Delta_1 \circ i_1=
(\overline{\sha}^\mathrm{q}(f) \otimes
\overline{\sha}^\mathrm{q}(f))\delta_1= (f \otimes f) \delta_1$.
Similarly, we find that $\epsilon_2 \circ
\overline{\sha}^\mathrm{q}(f) = \epsilon_1$ holds. Therefore,
$\overline{\sha}^\mathrm{q}(f)$ is a coalgebra morphism, which
concludes the proof.
\end{proof}

In the following we denote the algebra product in $H$ by $m_H(a
\otimes b) = a \cdot b$. Also, recall that $H \subset T^{+}(H)$
always means $i_H(H) := H \otimes 1_{\mathbb{K}} \subset
T^{+}(H)$. We omit the inclusion map $i$ for transparency.

Since $P_H$ is the right-shift operator, one can show the
following

\begin{prop} \label{prop:coideal}
For any words $h_1 \otimes \ldots \otimes h_n \in
(T^+(H),\bar{\bullet}^q,P_H)$, $q \in \{ {\mathrm{r}},\ell \}$ we
have
\begin{equation}
\label{rel:counit}
    \epsilon(h_1 \otimes \ldots \otimes h_n) = \epsilon(h_1) \cdots \epsilon(h_n)
\end{equation}
and for $q=\theta$ we find
 $\epsilon(h_1 \otimes \ldots
\otimes h_n) =  (-\theta)^{n-1} \epsilon(h_1) \cdots
\epsilon(h_n).$ For words $h_1 \otimes \ldots \otimes h_n \in
(T^+(H),\bar{\bullet}^q,P_H)$, $q \in \{ \theta, {\mathrm{r}},\ell
\}$
we have
\begin{equation}
\label{rel:coprod}
   \Delta(h_1 \otimes \ldots \otimes h_n)
             = (h_1)_{(1)} \ldots (h_n)_{(1)} \otimes \bigl((h_1)_{(2)}  \cdot \ldots  \cdot (h_n)_{(2)}\bigr),
\end{equation}
so that in fact $T^+(H)$ forms a right comodule over $H$, i.e.
$\Delta: T^+(H) \rightarrow T^+(H) \otimes H$.

For $h_1 \otimes \ldots \otimes h_n \in
(T^+(H),\bar{\bullet}^q,P_H)$, $q \in \{ {\mathrm{r}},\ell \}$
$$
    (\epsilon \otimes \id_{T^+(H)}) \Delta(h_1 \otimes \ldots \otimes h_n)
                           =  \bigl(h_1 \cdot \ldots \cdot  h_n \bigr),
$$
and for $q=\theta$ we find
 $   (\epsilon \otimes \id_{T^+(H)}) \Delta(h_1 \otimes \ldots \otimes h_n)
                           =  (-\theta)^{n-1}\bigl(h_1 \cdot \ldots \cdot  h_n
                           \bigr).$
\end{prop}

\begin{proof}
The proof follows by some algebra gymnastic. The last equation
comes from $(\mathbb{K}, -\theta \id_{\mathbb{K}})$ being an
$\textrm{{\bf{R}}}$-algebra of weight $\theta$.
\end{proof}

Moreover, as $T^+(H)$ is a bimodule over $H$, i.e., for $h \in H$
and $X \in T^+(H)$ we have
$$
    hX := i_H(h) \:\bar{\bullet}^q\: X \in T^+(H)
       \quad \; {\rm{ and }}\; \quad
    Xh := X \:\bar{\bullet}^q\: i_H(h) \in T^+(H),
$$
we observe that $(T^+(H),\Delta)$ is a right-covariant
bimodule~\cite{KS}, i.e.,
$$
    \Delta(hX + X'h') = (i_H^{\otimes 2} \circ \delta(h)) \:\bar{\bullet}^q\: \Delta(X) +
                        \Delta(X') \:\bar{\bullet}^q\: (i_H^{\otimes 2} \circ \delta(h')),
$$
and
$$
    (\Delta \otimes \id_{T^+(H)})\Delta = (\id_{T^+(H)} \otimes \delta)\Delta
   \quad \; {\rm{ and }}\; \quad
    (\id_{T^+(H)} \otimes \epsilon)\Delta = \id_{T^+(H)}.
$$
Defining instead a left $\textrm{{\bf{R}}}$-bialgebra on $
(T^+(H),\bar{\bullet}^q,P_H)$, we find a that $(T^+(H),\Delta)$ is
a left-covariant bimodule.\smallskip

\noindent Denote by $\mathcal{P}(C)$ the set of primitive elements
of the coalgebra $C$.

\begin{cor}
Let $H$ be a unital commutative bialgebra. Then,
$\mathcal{P}((T^+(H),\bar{\bullet}^q ,P_H)) = \mathcal{P}(H)$.
\end{cor}

\begin{proof}
This is due to the fact that $\Delta: T^+(H) \rightarrow T^+(H)
\otimes H$.
\end{proof}


\subsection{Case 2: unital algebras}
\label{ssect:case2}

Let $A$ be a commutative unital algebra with unit $1_A$. We would
like to establish a result similar to Theorem~\ref{thmHopfnij} for
commutative algebras. However, as the reader will observe, we must
restrict ourselves to the (right) Nijenhuis bialgebra framework.

We will use a somewhat unusual action of the algebra unit $1_A$
which will be explained further below. The following result was
found by one of us in \cite{Lerouxarth2}.

\begin{prop} \label{prop:newProd}
Let $(A,m_A, 1_A)$ and $(B,m_B,1_B)$ be two unital associative
algebras. Then, $A \circledcirc B$ whose product $\amalg$ is
defined as follows,
 \allowdisplaybreaks{
\begin{eqnarray*}
    a_1 \otimes b_1 \amalg a_2 \otimes 1_B &=& 0 \quad \textrm{if $b_1 \neq 1_B$ and $a_2 \neq 1_A$,}\\
    a_1 \otimes 1_B \amalg a_2 \otimes b_2 &=& 0 \quad \textrm{if $a_1 \neq 1_A$ and $b_2 \neq 1_B$, }\\
    a_1 \otimes b_1 \amalg a_2 \otimes b_2 &=& m_A(a_1, a_2) \otimes m_B(b_1,b_2) \qquad \textrm{otherwise,}
\end{eqnarray*}}
is associative and has $1_A \otimes 1_B$ as unit. It is isomorphic
to $A \otimes B$ as $\mathbb{K}$-vector space and, moreover, if
$A$ and $B$ are commutative so is the product $\amalg$.
\end{prop}

\begin{proof} The proof of the associativity follows easily by
direct calculations of the different cases. Commutativity is
readily verified.
\end{proof}

\begin{cor} \label{cor:newProd}
Let
$\overline{\sha}^\mathrm{r}(A)=(T^+(A),\bar{\bullet}^{\mathrm{r}}
,P_A)$ be the free commutative Nijenhuis algebra over a
commutative algebra $A$ with unit $1_A$. Then, the following
extension of $\bar{\bullet}^{\mathrm{r}}$ to
$\overline{\sha}^\mathrm{r}(A)^{\otimes 2}$, still denoted by
$\bar{\bullet}^{\mathrm{r}}$ and defined by
 \allowdisplaybreaks{
\begin{eqnarray*}
    a_1 \otimes b_1 \:\bar{\bullet}^{\mathrm{r}}\: a_2 \otimes 1_A &=&
                        0 \quad \textrm{if $b_1$ and $a_2$ are both different from $1_A$, }\\
    a_1 \otimes 1_A \:\bar{\bullet}^{\mathrm{r}}\: a_2 \otimes b_2 &=&
                        0 \quad \textrm{if $a_1$ and $b_2$ are both different from $1_A$, }\\
    a_1 \otimes b_1 \:\bar{\bullet}^{\mathrm{r}}\: a_2 \otimes b_2 &=&
                        a_1 \bar{\bullet}^{\mathrm{r}} a_2 \otimes b_1 \bar{\bullet}^{\mathrm{r}} b_2 \quad \textrm{otherwise,}
\end{eqnarray*}}
is commutative, associative, and has $1_A \otimes 1_A$ as unit.
\end{cor}

\begin{prop} \label{newNijenhuis}
Let $A$ be a (commutative) unital algebra with unit $1_A$. Then,
$\overline{\sha}^\mathrm{r}(A)^{\otimes 2}$ with the product as in
Corollary~\ref{cor:newProd} is a Nijenhuis algebra when equipped
with the Nijenhuis map $\mathsf{P}_A := P_A \otimes \id_A$.
\end{prop}

\begin{proof} First we should remark, that the statement of this
proposition extends naturally to
$\overline{\sha}^\mathrm{r}(A)^{\otimes n}$ for $n>2$.

In the following we reserve the tensor product sign for elements
in $T^+(A)^{\otimes 2}$ and use the word notation for elements in
$T^+(A)$. We have to prove that for $a,b,c,d \in T^+(A)$
$$
        \mathsf{P}_A(a \otimes b) \bar{\bullet}^{\mathrm{r}} \mathsf{P}_A(c \otimes d)
                            + \mathsf{P}_A^2 \bigl((a \otimes b) \bar{\bullet}^{\mathrm{r}} (c \otimes d)\bigr)
        = \mathsf{P}_A\bigl(\mathsf{P}_A(a \otimes b) \bar{\bullet}^{\mathrm{r}} (c \otimes d)
                            + (a \otimes b) \bar{\bullet}^{\mathrm{r}} \mathsf{P}_A(c \otimes d) \bigr),
$$
that is,
$$
    (1_A a \otimes b) \bar{\bullet}^{\mathrm{r}} (1_A c \otimes d)
                          + \mathsf{P}_A^2\bigl( (a \otimes b) \bar{\bullet}^{\mathrm{r}} (c \otimes d) \bigr)
    = \mathsf{P}_A\bigl((1_A a \otimes b)\bar{\bullet}^{\mathrm{r}} (c \otimes d)
                          + (a \otimes b) \bar{\bullet}^{\mathrm{r}} (1_A c \otimes d)\bigr).
$$
This equality obviously holds if $a,b,c$ and $d$ are different
from $1_A$. We are left to check the two cases.
\begin{itemize}
    \item {\textsf{Case} $d = 1_A$ and $b \not= 1_A$.
         If $c \not = 1_A$, the equality turns out to be true.
         If $c = 1_A$ then we find, $0 + P_A^2(a) \otimes b = 0 + P_A^2(a) \otimes b.$}

    \item {\textsf{Case} $b = 1_A$ and $d \not= 1_A$. This case follows by an argument as above.}

    \item {\textsf{Case} $b=1_A$ and $d=1_A$. This case is readily verified.}
\end{itemize}
Consequently, $\mathsf{P}_A$ is a Nijenhuis operator on
$\overline{\sha}^\mathrm{r}(A)^{\otimes 2}$ with the product as in
Corollary~\ref{cor:newProd} .
\end{proof}

\begin{thm} \label{thm:conij}
Let $A$ be a commutative unital algebra with unit $1_A$. We denote
the algebra product by brackets, $[-;-]_A$. We assume that no $a
\not=1_A \in A$ has an inverse (e.g. take the free commutative
algebra over a $\mathbb{K}$-vector space $V$). Then, the free
commutative Nijenhuis algebra
$(T^+(A),\bar{\bullet}^{\mathrm{r}},P_A)$ can be turned into a
right Nijenhuis bialgebra.
\end{thm}

\begin{proof}
We found in Proposition~\ref{newNijenhuis} that $(T^+(A)^{\otimes
2},\bar{\bullet}^{\mathrm{r}},\mathsf{P}_A)$ is also a Nijenhuis
algebra. Consider the following commutative diagram
\begin{diagram}
    A          & \rInto^{\phantom{mim} i_A \phantom{mim}} & (T^+(A),\bar{\bullet}^{\mathrm{r}},P_A)\\
               & \rdTo_{\delta\ }                         & \dTo_{\Delta }\\
               &                                          & (T^+(A)^{\otimes 2},\bar{\bullet}^{\mathrm{r}},\mathsf{P}_A)
\end{diagram}

\noindent where the inclusion is defined as
before~(\ref{inclusionmap}), i.e., $i_A(a):= a \otimes
1_{\mathbb{K}} \in T^+(A)$, giving an algebra morphism. The map
$\delta$ is defined by $a \mapsto i_A(1_A) \otimes i_A(a) + i_A(a)
\otimes i_A(1_A)$ if $a \not=1_A$ and $1_A \mapsto i_A(1_A)
\otimes i_A(1_A)$ otherwise. We will omit the inclusion map $i_A$
in the following. The map $\delta$ so defined is a unital algebra
morphism since $\delta(a) \bar{\bullet}^{\mathrm{r}} \delta(b) = a
\bar{\bullet}^{\mathrm{r}} b \otimes 1_A + 1_A \otimes a
\bar{\bullet}^{\mathrm{r}} b = [a;b]_A \otimes 1_A + 1_A \otimes
[a;b]_A = \delta([a;b]_A)$. Two remarks are in order. First, the
particular product in Proposition~\ref{prop:newProd} applies at
this point. Secondly, the special assumption that no element in
$A$ has an inverse is needed here for consistency.

Consequently, there exists a unique Nijenhuis algebra map $\Delta$
which extends $\delta$, such that $\Delta \circ i_A = \delta$.
Observe that our construction extends to $(T^+(A)^{\otimes
3},\bar{\bullet}^{\mathrm{r}},\mathsf{P}_A \otimes \id_A)$ as
well, making it into a commutative Nijenhuis algebra. As before
let us define the following map $\delta^{(2)}:= i_A^{\otimes 3}
\circ (\id_A \otimes \delta)\delta: A \to T^+(A)^{\otimes 3},$
which is a unital algebra morphism and we have $\delta^{(2)}=
i_A^{\otimes 3} \circ (\delta \otimes
\id_A)\delta=:\hat{\delta}^{(2)}$. Define
$\Delta^{(2)}:=(\id_{T^+(A)} \otimes \Delta)\Delta$ and
$\hat{\Delta}^{(2)}:=(\Delta \otimes \id_{T^+(A)})\Delta$ and
consider the commutative diagram
\begin{diagram}
    A          & \rInto^{\phantom{mim} i_A \phantom{mim}} & (T^+(A),\bar{\bullet}^{\mathrm{r}},P_A)\\
               & \rdTo_{\delta^{(2)}\ }                   & \dTo_{\Delta^{(2)}}\\
               &                                          & (T^+(A)^{\otimes 2},\bar{\bullet}^{\mathrm{r}},\mathsf{P}_A)
\end{diagram}

\noindent In fact, the map $\delta^{(2)}(a) =
\hat{\delta}^{(2)}(a) = a \otimes 1_A \otimes 1_A + 1_A \otimes a
\otimes 1_A + 1_A \otimes 1_A\otimes a $ for $a \not = 1_A$ and
$1_A \mapsto 1_A \otimes 1_A \otimes 1_A$ is an associative unital
algebra morphism which is extended by $(\Delta \otimes
\id_{T^+(A)}) \Delta$ (as well as by $(\id_{T^+(A)} \otimes
\Delta) \Delta$), hence the coassociativity of the coproduct
$\Delta$ follows by the unicity of the extension. Consequently,
the map $\Delta$ is both a coassociative coproduct and a Nijenhuis
algebra morphism.

Consider the unital algebra map $\epsilon: A \rightarrow
\mathbb{K}$ defined by $\epsilon(1_A)=1_\mathbb{K}$ and zero
otherwise. We get as before by extension a counit, still denoted
by $\epsilon$, satisfying the equations
$$
    (\id_{T^+(A)} \otimes \epsilon)\Delta= \id_{T^+(A)}, \quad \textrm{and} \quad
    (\epsilon \otimes \id_{T^+(A)})\Delta =[-]_A^*,
$$
where $[-]_A^*:T^+(A) \rightarrow i(A)$ is inductively defined by,
$[a_1]_A^*:=a_1 \otimes 1_{\mathbb{K}}$ and $[a_1 \otimes \dots
\otimes a_n]^*_A:=[[a_1 \otimes \dots \otimes a_{n-1}]_A^*; a_n]_A
\otimes 1_{\mathbb{K}}$.
\end{proof}

\noindent We hence find that
$$
    \Delta(x) = x \otimes 1_A + 1_A^{\otimes l(x)} \otimes [x]^*_A, \qquad\
    \Delta(1_A^{\otimes n})= 1_A^{\otimes n} \otimes 1_A,
$$
where $l(x)=l(x_1 \otimes \dots \otimes x_n):=n$ is the length of
$x \in T^+(A)$.

\noindent Denote by $\mathcal{P}(C)$ the set of primitive elements
of the coalgebra $C$.

\begin{cor}
Let $A$ be a unital commutative algebra. Then
$\mathcal{P}((T^+(A),\bar{\bullet}^{\mathrm{r}} , P_A))=A$.
\end{cor}

\vspace{1cm}

\noindent \textbf{Acknowledgments:} The first author acknowledges
greatly the support by the European Post-Doctoral Institute and
Institut des Hautes \'Etudes Scientifiques (I.H.\'E.S.). L.~P.
would like to acknowledge the warm hospitality he experienced
during his stay at I.H.\'E.S. when major parts of this paper were
written. We thank W.~Schmitt and D.~Manchon for useful discussions
and comments.


\end{document}